\pdfoutput=1

\documentclass[letterpaper, 12pt]{article}

\usepackage{amssymb,amsmath,amsthm}
\usepackage{epsf}
\usepackage{times,bm,mathrsfs}
\usepackage{amsmath}
\usepackage{amssymb}
\usepackage{amsfonts}
\usepackage{mathrsfs}
\usepackage{bbm}
\usepackage{color}
\usepackage{graphicx}
\usepackage{amscd}
\usepackage{epsfig}

\usepackage{mathrsfs}  
\usepackage{color}
\usepackage[T1]{fontenc}
\usepackage[utf8]{inputenc}
\usepackage[english]{babel}

\usepackage{graphicx}
\usepackage{amssymb}
\usepackage{mathrsfs}
\usepackage{wrapfig}

\usepackage{booktabs}

\definecolor{Red}{rgb}{1,0,0}
\definecolor{Blue}{rgb}{0,0,1}
\definecolor{Olive}{rgb}{0.41,0.55,0.13}
\definecolor{Green}{rgb}{0,1,0}
\definecolor{MGreen}{rgb}{0,0.8,0}
\definecolor{DGreen}{rgb}{0,0.55,0}
\definecolor{Yellow}{rgb}{1,1,0}
\definecolor{Cyan}{rgb}{0,1,1}
\definecolor{Magenta}{rgb}{1,0,1}
\definecolor{Orange}{rgb}{1,.5,0}
\definecolor{Violet}{rgb}{.5,0,.5}
\definecolor{Purple}{rgb}{.75,0,.25}
\definecolor{Brown}{rgb}{.75,.5,.25}
\definecolor{Grey}{rgb}{.5,.5,.5}
\definecolor{Black}{rgb}{0,0,0}

\newcommand{\R}{{\mathbb{R}}}
\newcommand{\E}{{\mathbb{E}}}
\newcommand{\pr}{{\mathbb{P}}}
\newcommand{\G}{{\mathbb{G}}}
\newcommand{\1}{{\textbf{1}}}

\newcommand{\T}{{\mathbb{T}}}
\newcommand{\tp}{$(\T,P)\textit{-walk on } G$ }
\newcommand{\tpns}{$(\T,P)\textit{-walk on } G$}

\newcommand{\gn}{N}

\newcommand{\D}{{\mathscr{D}}}

\newcommand{\A}{{\mathscr{A}}}
\newcommand{\sP}{{\mathscr{P}}}

\newtheorem{theorem}{Theorem}[section]
\newtheorem{definition}{Definition}

\newtheorem{prop}{Proposition}[section]
\newtheorem{remark}{Remark}[section]

\newcommand{\mutrue}{\mu_{\mathrm{true}}}
\renewcommand{\deg}{\mathrm{deg}}
\newcommand{\var}{\mathrm{Var}}
\newcommand{\cov}{\mathrm{Cov}}
\newcommand{\rds}{\mathrm{RDS}}
\newcommand{\gls}{\mathrm{GLS}}
\newcommand{\rank}{\mathrm{rank}}
\newcommand{\auto}{\mathrm{auto}}
\newcommand{\diag}{\mathrm{diag}}
\newcommand{\sbm}{\mathrm{SBM}}
\newcommand{\true}{\mathrm{true}}
\newcommand{\vh}{\mathrm{VH}}
\renewcommand{\ss}{\mathrm{SS}}
\newcommand{\rse}{\mathrm{RSE}}
\newcommand{\np}{\mathrm{NP}}
\renewcommand{\cov}{\mathrm{Cov}}

\addtocounter{page}{-1}

\begin{document}


\title{Generalized least squares can overcome\\ the critical threshold\\ in respondent-driven sampling
}

\author{
	Sebastien Roch\footnote{Department of Mathematics, UW--Madison.}
		\and
		Karl Rohe\footnote{Department of Statistics, UW--Madison. To whom correspondence should be addressed. E-mail: karlrohe at stat dot wisc dot edu}
}

\maketitle

\begin{abstract}
	In order to sample marginalized and/or hard-to-reach populations, resp\-on\-dent-driven sampling (RDS) and similar techniques reach their participants via peer referral.
	Under a Markov model for RDS, previous research has shown that if the typical participant refers too many contacts, then the variance of common estimators does not decay like $O(n^{-1})$, where $n$ is the sample size. This implies that confidence intervals will be far wider than under a typical sampling design.  Here we show that generalized least squares (GLS) can effectively reduce the variance of RDS estimates.  In particular, a theoretical analysis indicates that the variance of the GLS estimator is $O(n^{-1})$.  We then derive two classes of feasible GLS estimators.  The first class is based upon a Degree Corrected Stochastic Blockmodel for the underlying social network. The second class is based upon a rank-two model.  It might be of independent interest that in both model classes, the theoretical results show that it is possible to estimate the spectral properties of the population network from the sampled observations.  Simulations on empirical social networks show that the feasible GLS (fGLS) estimators can have drastically smaller error and rarely increase the error.  A diagnostic plot helps to identify where fGLS will aid estimation.  The fGLS estimators continue to outperform standard estimators even when they are built from a misspecified model and when there is preferential recruitment.  
\end{abstract}

\thispagestyle{empty}

\clearpage

\section{Introduction}

Respondent-driven sampling (RDS) is a popular network-based approach to sample marginalized and/or hard-to-reach populations \cite{heckathorn1997respondent}. 
RDS has become particularly popular in HIV research because the populations most at risk for HIV (e.g., people who inject drugs, female sex workers, and men who have sex with men) cannot be sampled using conventional techniques.  Several  domestic and international institutions use RDS to quantify the prevalence of HIV in at-risk populations, including the Centers for Disease Control (CDC), the World Health Organization (WHO), and the Joint United Nations Programme on HIV/AIDS (UNAIDS) \cite{WHO}.   The most recent review of the literature in 2015 counted over 460 different RDS studies, in 69 different countries \cite{white2015strengthening}.

Because RDS collects samples from link-tracing the relationships in a social network, adjacent samples are dependent.  In a simulation study, \cite{goel2010assessing} showed how this can lead to highly variable estimates.   
Under independent sampling, the variance of standard estimators decays like $O(n^{-1})$.  This implies that a sample size of $4n$ will have a 50\% smaller standard error than a sample of size $n$.  However, this does not necessarily hold for RDS.  
Under a Markov model, \cite{treevar} 
showed how the dependence induced by RDS can drastically inflate the variance of traditional estimators, making it decay at a rate slower than $O(n^{-1})$.  
This implies that reducing the sampling error by 50\% can require far more than 4 times as many samples. This means that confidence intervals are much wider than under independent sampling. 
Using the covariance function derived in \cite{treevar}, this paper studies the generalized least squares (GLS) estimator for RDS.  Our theoretical analysis establishes that the variance of the GLS estimator is $O(n^{-1})$.  
We then derive a feasible GLS (fGLS) estimator based upon the Degree Corrected Stochastic Blockmodel.  Two alternative estimators are derived.  These estimators first construct estimates about the spectral properties of the population social graph, which might be of independent interest.
Our fGLS estimators easily accommodate any preliminary re-weighting of the data to adjust for the sampling biases that occur in RDS (e.g. \cite{volz2008probability, gile2011improved}). We study these estimators with simulations and propose a simple diagnostic plot to compare the different fGLS estimators.

\begin{figure}[htbp] 
   \centering
   \includegraphics[height=2.8in]{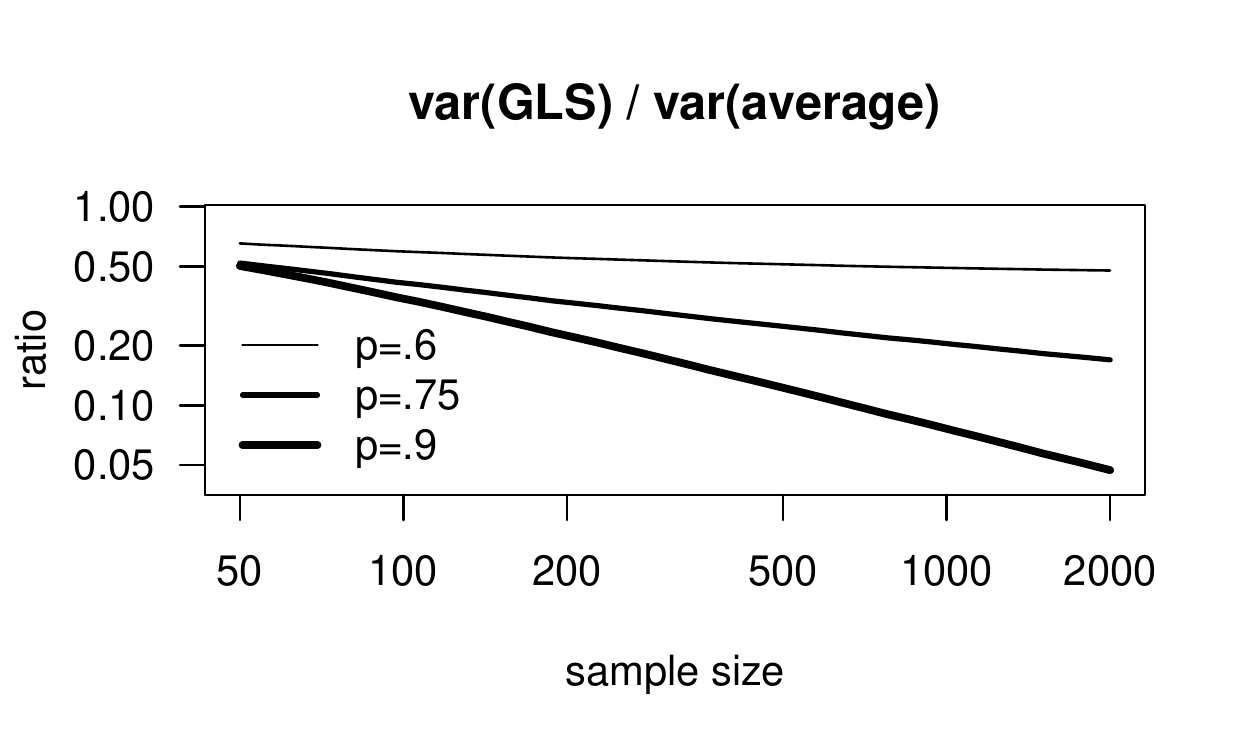} 
   \caption{
In this experiment, GLS provides dramatic improvements when the sample is large and the correlation between samples (i.e. $p$) is high. 
  Both axes are on the log scale.
}
   \label{fig:motivate}
\end{figure}

\paragraph{A simple motivating example}  Figure \ref{fig:motivate}  
uses a  model studied in \cite{goel2009respondent}.  In this example, the population that we wish to sample is equally divided into two groups: HIV+ and HIV-.  The seed participant is selected uniformly at random. Starting with the seed participant, every participant refers two additional participants (as in a complete binary tree).  
 The participant refers a person that matches their own HIV status with probability $p$ and refers a person with the opposite status with probability $1-p$.  Each referral is independent.
Using this sample, we wish to estimate the proportion of the population that is HIV+ (i.e., 50\%).  Figure \ref{fig:motivate} compares two estimators, (i) the sample proportion and (ii) the GLS estimator proposed in this paper.

Under this sampling with replacement model, the variance of both the sample proportion and the GLS estimator have closed form solutions. (See \cite{goel2009respondent} and the proof of Theorem \ref{thm:oneEigen}.)
Figure \ref{fig:motivate} gives the ratio of these formulas as a function of the sample size $n$.
There are three lines,  corresponding to $p= .6, p = .75,$ and $p=.9$. In all cases, the lines are less than one, indicating that the GLS estimator has a smaller variance than the sample proportion. 
Under this simulation model, if $p>.86$, then the variance of the sample proportion decays slower than $O(n^{-1}$) \cite{goel2009respondent, treevar}.
As Theorems \ref{thm:main} and \ref{thm:oneEigen} below show,  the variance of the GLS estimator converges to zero like $O(n^{-1})$.  So, as $n$ increases, the bottom line converges to zero. The other two lines, on the other hand, do not converge to zero.

\section{Preliminaries} \label{sec:prelim}

The Markov model used in this paper is a straightforward combination of the Markov models developed in the RDS literature (e.g. \cite{heckathorn1997respondent, salganik2004sampling, volz2008probability} and \cite{goel2009respondent}). 

The social network, $G = (V,E)$ consists of the node set $V = \{1, \dots, N\}$ and the edge set $E = \{(i,j) : \mbox{ $i$ and $j$ can refer each other}\}$. To simplify notation, $i \in G$ is used synonymously  with $i \in V$.  Unless otherwise noted, everything below also applies to weighted graphs.   Let $w_{ij}$ be the weight of edge $(i,j) \in E$, which models preferential recruitment as described in Section~\ref{sec:pref}. If $(i,j) \not \in E$, define $w_{ij} = 0$. If the graph is unweighted, then let $w_{ij} = 1$ for all $(i,j) \in E$.  Throughout this paper, the graph is undirected, i.e., $w_{ij} = w_{ji}$ for all pairs $i,j$. Define the degree of node $i$ as 
$\deg(i) = \sum_j w_{ij}$.
For each node $i \in G$, let $y(i) \in \R$ denote some characteristic of this node (e.g., the indicator of HIV status).  We wish to estimate the population average 
\begin{equation}
\label{eq:mutrue}
\mutrue = \frac{1}{\gn} \sum_{i \in G}  y(i).
\end{equation}

We assume that the nodes are sampled with a Markov process that is indexed by a rooted tree $\T$ (i.e., a connected graph with $n$ nodes, no cycles, and a vertex $0$).  The seed participant is vertex $0$ in $\T$. 
To simplify notation, $\sigma \in \T$ is used synonymously with $\sigma$ belonging to the vertex set of $\T$.  For any node in the tree $\sigma \in \T$,  denote $\sigma' \in \T$ as the parent of $\sigma$ (the node one step closer to the root). 
 Define the matrix $P \in \R^{\gn \times \gn}$ as
\begin{equation}\label{eq:pdef}
P_{ij} = \frac{w_{ij}}{\deg(i)}.
\end{equation}
Because the graph is undirected, $P$ is a reversible Markov transition matrix with a stationary distribution $\pi: G \rightarrow \R$.  
Our sample is the set of random nodes 
 \[\{X_\sigma \in G: \sigma \in \T\},\]
 where $X_0$ is initialized with $\pi$ and each transitions $X_{\sigma'} \rightarrow X_{\sigma}$ is independent with
\[\pr(X_\sigma = j | X_{\sigma'} = i) = P_{ij}, \mbox{ for } i,j \in G.\]
Observe that $\T$ and $G$ are distinct graphs: the nodes in $\T$ index the Markov process while the nodes in $G$ are its state space.
 Following \cite{benjamini1994markov}, we refer to this stochastic process as a \tpns. 

When the $X_\tau$'s sample the target population in $G$, we observe
\[Y_\tau = y(X_\tau) \ \mbox{ for } \tau \in \T.\] 
Under the stationary \tpns, the sample average of the $Y_\tau$'s is an estimate of 
\[\mu = \E(Y_0) = \sum_i y(i) \pi_i.\]  
In general, $\mu \ne \mutrue$ (where $\mutrue$ was defined~\eqref{eq:mutrue}) and the sample average must be adjusted with sampling weights in order to obtain an unbiased estimator of $\mutrue$.  
Define 
\[y^\pi(i) = \frac{y(i)}{\pi_{i}N}, \quad\quad Y_\tau^\pi = y^\pi(X_\tau).\]  
The sample average of the $Y_\tau^\pi$'s is the inverse probability weighted (IPW)  estimator; it is an unbiased estimator of $\mutrue$ \cite{horvitz1952generalization}.  However, the weights $\pi_i N$ are unknown and must be estimated with additional information, as we describe next. 

Under the \tpns, 
\[N \pi_i  = \frac{N \deg(i)}{\sum_j \deg(j)} = \frac{N\deg(i)} {N \bar d} = \frac{\deg(i)}{ \bar d},\] 
where $\bar d =  N^{-1} \sum_j \deg(j)$.
The popular Volz-Heckathorn estimator replaces $\bar d$ with the harmonic mean of those degrees \cite{volz2008probability}. Recall that $\T$ has $n$ nodes and define 
\[H_{\deg}^{-1} = \frac{1}{n} \sum_{\tau \in \T} \frac{1} {\deg(X_\tau)}, \  \ 
\hat \pi_i = H_{\deg}^{-1}\deg(i), \  \
y^{\hat \pi}(i) = \frac{y(i)}{\hat \pi_{i}},
\] 
and $Y_\tau^{\hat \pi} = y^{\hat \pi}(X_\tau)$.
The Volz-Heckathorn estimator is the sample average of the $Y_\tau^{\hat \pi}$'s and it is an asymptotically unbiased estimator of $\mutrue$ under the {\it $(\T,P)$-walk on $G$}.\footnote{In practice, $\deg(i)$ is estimated by asking participants how many contacts they have.  Recall that $\deg(i) = \sum_j w_{ij}$.  If the graph is weighted, 
then the \tp exhibits preferential recruitment (as discussed in Section S1) and the number of contacts will not necessarily align with $\deg$, making the estimator biased.}

\begin{remark} \label{remark:unbiasWeighting} The next section will drop the superscript $\pi$ and $\hat \pi$ in $Y_\tau^{\pi}$ and $Y_\tau^{\hat \pi}$.  Using the $Y_\tau$'s to construct the GLS estimator will lead to an unbiased estimator of $\mu$.  In practice, before doing any of the GLS computations, one could replace the $Y_\tau$'s with $Y_\tau^{\pi}$ or $Y_\tau^{\hat \pi}$ in order to estimate $\mutrue$.   The simulations in this paper use a reweighting that is similar to $\hat \pi$, but replaces $H_{\deg}^{-1}$ with a GLS estimate of $\E(1/\deg(X_\tau))$.  
In~\cite{gile2011improved}, sampling weights are estimated under an alternative, non-Markovian model.  These weights could also be used before doing GLS computations.
\end{remark}

\section{GLS for RDS}

The GLS estimator is the weighted average of the $Y_\tau$'s with smallest variance \cite{aitken1936iv}, i.e., it is the solution $g^*$ to
\begin{equation}\label{eq:glsobjective}
\min_{g} \var\left(\sum_{\tau \in \T} g_\tau Y_\tau\right) \quad \mbox{ such that } \quad \sum_{\tau \in \T} g_\tau = 1.
\end{equation}
Because of the constraint that the weights $g_\tau$ sum to one, the linearity of expectation, and the fact that the \tp is stationary, the resulting estimator is an unbiased estimate of $\E(Y_\tau)$.
Define the covariance matrix $\Sigma \in \R^{n \times n}$ as
\begin{equation}\label{eq:covmat}
\Sigma_{\sigma,\tau} = \cov_{\rds}(Y_\sigma, Y_{\tau}),
\end{equation}
which is assumed to be nonsingular.
It can be seen that the solution to \eqref{eq:glsobjective} depends upon solving a system of equations involving the covariance matrix, namely that  $g^* = (x^T \1)^{-1} x^T$ where $\Sigma x = \1$.  
(Throughout, we use the notation $\1_M$ for the all-one vector of length $M$. We drop the length when clear from context.)
If $Y \in R^n$ is the vector of $Y_\tau$'s, then the GLS estimator can be expressed as
\begin{equation}\label{eq:gls}
\hat \mu_{\gls} = (\1^T \Sigma^{-1} \1)^{-1} \1^T \Sigma^{-1} Y.
\end{equation}

\vspace{-.01 in}


The rest of this section contains our main theoretical results, which study how
\begin{equation} \label{eq:glsvar}
\var_{\rds}(\hat \mu_{\gls})   = (\1^T \Sigma^{-1} \1)^{-1}
\end{equation}
 decays with the sample size.

\subsection{Main result}   \label{sec:main}

In our main result, we assume that $\T$ is a complete
binary tree with $n$ nodes, but we expect the result to hold for more general tree topologies.
\begin{theorem}[Main Result] \label{thm:main}
Let $\{X_\tau : \tau \in \T\}$ be sampled from the $(\T,P)$\textit{-walk on} $G$ for a fixed $N\times N$ transition matrix $P$ that is irreducible and reversible with respect to a stationary distribution $\pi$.  If  $\T$ is a complete binary tree with $n$ nodes, then the variance of the GLS
estimator defined in \eqref{eq:gls} decays like $O(n^{-1})$
as $n \rightarrow \infty$.  
\end{theorem}
The proof, which is contained in Section~\ref{sec:technical}, does not directly compute the variance of the GLS estimator. Instead, it proceeds by constructing an explicit linear estimator and relies on the variational characterization \eqref{eq:glsobjective} of $\hat \mu_{\gls}$. We emphasize that computing $\hat \mu_{\gls}$ requires the covariance matrix $\Sigma$, which is typically unknown.  The next section proposes a technique to estimate $\Sigma$ that is based upon the Stochastic Blockmodel. We also point out that the result in Theorem \ref{thm:main} is asymptotic and, as such, is only meaningful for $n$ large enough.

Before moving on to practical estimators, we give a more precise result on the constant in the $O(n^{-1})$ by making further assumptions on the spectral properties of $P$ or of the features $y$.  
The eigenvectors of the reversible transition matrix $P$, denoted $f_1, \dots, f_{\gn}: V \rightarrow \R$, are real-valued functions of the nodes $i \in G$ that are orthonormal with respect to the inner product 
\begin{equation} \label{def:inner}
\langle f_a, f_b \rangle_\pi = \sum_{i \in G} f_a(i) f_b(i) \pi_i.
\end{equation}
(See, e.g., Lemma 12.2 of \cite{levin2009markov}.)
We take the eigenfunction $f_1$ corresponding to the eigenvalue $1$ to be the constant vector $\1$.
Define $\beta_\ell = \langle y, f_\ell \rangle_\pi$ for $\ell = 1, \dots, N$
and note that $\mu = \beta_1 = \sum_i y(i)\pi_i$.  Let $\lambda_1, \dots, \lambda_\gn$ be the eigenvalues of $P$ corresponding to $f_1, \dots, f_\gn$.  For each node $i \in G$,  $y$ decomposes as follows
\begin{equation}\label{eq:y-eigendecomposition}
y(i) 
= \mu +\sum_{\ell =2}^{N} \langle y, f_\ell \rangle_\pi 
f_\ell(i) 
= \sum_{\ell =1}^{N} \beta_\ell
f_\ell(i).
\end{equation}
Under the \tpns, the covariance is stationary with  auto-covariance function
\begin{equation}\label{eq:gammadef}
\gamma( d) = \sum_{\ell=2}^{\gn} \beta_\ell^2 \lambda_\ell^{d}.  
\end{equation}
That is, the covariance matrix has the form $\Sigma_{\sigma,\tau} = \gamma(d(\sigma, \tau))$, 
%
%
%
where  $d(\sigma, \tau)$ is the graph distance (i.e. minimum path length) between $\sigma$ and $\tau$ in $\T$ \cite{treevar}.

When the auto-covariance further simplifies to
\begin{equation}\label{eq:rank2model}
\gamma(d) =  \beta^2 \lambda^{d},
\end{equation}
for some $\lambda, \beta \in \R$, then we call the \tp with feature $y$ a \textbf{rank-two model}.  
%
For instance, if $\rank(P) = 2$, then as the name suggests,  we have a rank-two model. In particular, 
all of the results in \cite{goel2009respondent} are for such transition matrices. 
Figure \ref{fig:motivate} in the introduction also studies such a rank-two model on two groups of people.  There are other sufficient conditions for \eqref{eq:rank2model}. For instance, if $y(i) = \mu + \beta_\ell f_\ell(i)$ for all nodes $i \in G$, then we have a rank-two model because $\beta_j = \langle y, f_j\rangle_\pi = 0$ for $j \not \in  \{1, \ell\}$.

%

\begin{theorem} \label{thm:oneEigen}
Under a rank-two model, 
\begin{equation} \label{eq:glsvar1}
n \var(\hat \mu_{\gls}) 
\rightarrow \left(\frac{1 + \lambda}{1 - \lambda} \right) \beta^2 \mbox{ as $n \rightarrow \infty$}.
\end{equation}
\end{theorem}
This proof follows from the fact that under a rank-two model, $\Sigma^{-1}$ has a closed form expression  (see Section~\ref{sec:technical}).

\section{Using RDS to estimate the spectral properties of the graph for feasible GLS} \label{sec:fgls}
The feasible GLS (fGLS)  estimator depends upon an estimated covariance matrix $\hat \Sigma$ (e.g. see \cite{amemiya1985advanced}), 
\begin{equation}\label{eq:muhatSigma}
\hat \mu_{\gls}(\hat \Sigma) = (\1 \hat \Sigma^{-1} \1)^{-1} \1 \hat \Sigma^{-1} Y.
\end{equation}
With this notation, observe that $\hat \mu_{\gls} = \hat \mu_{\gls}(\Sigma)$.

In our setting, estimating $\Sigma$ is equivalent to estimating
$\gamma( \cdot)$.
We propose and compare several estimators for $\gamma$.
An estimator based upon the Degree Corrected Stochastic Blockmodel (DC-SBM) is derived in this section. Two additional estimators
based upon the rank-two assumption are derived in Section~\ref{sec:rank2est}.  The first rank-two estimator, $\hat \mu_{\auto}$, relies upon a plug-in estimator for the correlation between $Y_{\sigma'}$ and $Y_\sigma$ (i.e., the auto correlation at lag 1).  The second rank-two estimator, $\hat \mu_{\Delta}$, relies upon plug-in estimators for the first and second differences, $\E(Y_{\sigma'} - Y_\sigma)^2$ and  $\E(Y_{(\sigma')'} - Y_\sigma)^2$.  


\subsection{Estimating the spectral properties of a Stochastic Blockmodel from an RDS}\label{sec:sbmfgls}

The Degree Corrected Stochastic Blockmodel is a generalization of the Stochastic Blockmodel~\cite{holland1983stochastic, karrer2011stochastic}.  Both are models for a random network with community structure. As the name suggests, the degree corrected model allows for degree heterogeneity within the blocks.
\begin{definition}[Degree Corrected Stochastic Blockmodel] 
Partition the $N$ nodes into $K$ blocks with $z$ : $\{1, 2,..., N\} \rightarrow \{1, 2,..., K\}$ and assign each node $i$ a value $\theta_i > 0$ such that the $\theta$s sum to one within each block, i.e.,
\begin{equation} \label{eq:theta}
\sum_{i:z(i) = u} \theta_i  = 1, \quad \mbox{ for all $u \in \{1, \dots, K\}$.}
\end{equation}
The block membership of node $i$ is $z(i)$ and the parameter $\theta_i$ controls the degree heterogeneity within each block. Let ${\bf B}$ be a symmetric $K \times K$ matrix such that ${\bf B}_{ab} \ge 0$ for all $a,b \in 1, \dots, K$.
Under the Degree Corrected Stochastic Blockmodel (DC-SBM), 
\[\pr(\{i,j\} \in E) = \theta_i \theta_j {\bf B}_{z(i), z(j) }\] 
for all pairs $i,j = 1,2,...,N$ and each possible edge is independent. 
\end{definition}

In much of the previous literature on the DC-SBM, the full network is observed and we wish to estimate the partition  $z$.  In this paper, we presume that $z$ is observed on the sampled nodes in the \tp and we wish to estimate the spectral properties of $P$.  
This is reasonable in RDS because each participant takes a survey which records several salient demographic variables (e.g., gender, race, neighborhood, etc). In practice, the block labels should be chosen such that they are highly autocorrelated from one referral to the next.
Many RDS papers already report such statistics.  For example, 
the original RDS paper \cite{heckathorn1997respondent} presents four empirical transition matrices on four different demographic partitions (i.e., race, gender, drug preference, and location).

The derivations below condition on the block labels $z$; only the graph $G$ is random. Let $A \in \{0,1\}^{N\times N}$ be the (random) adjacency matrix; $A_{ij} = 1$ if and only if $(i,j) \in E$.  
Define $\A \in [0,1]^{N\times N}$ such that $\A_{ij} = \E(A_{ij}) = \pr( (i,j) \in E)$.   
  Define $\D \in \R^{N \times N}$ as a diagonal matrix with $(i,i)$-th element  $\sum_j \A_{ij}$.  Define $\sP = \D^{-1} \A$ as a population version of $P$.

The inspiration for the following estimators is based on a population version of the chain and relies on three results. 
Define the matrix $\hat Q \in \R^{K \times K}$ such that for any two blocks $u,v$,
\begin{equation}\label{Qdef}
\hat Q_{uv} = \frac{1}{n} \times \mbox{ number of referrals from block $u$ to block $v$}.
\end{equation}
Proposition \ref{prop:Q} below shows that $\hat Q$ is an estimator of $\bf B$ under a $(\T,\sP)\textit{-walk on } G$. Then, Proposition \ref{prop:B} shows that a normalized version of $\bf B$ has spectral properties that match the spectral properties of $\sP$.  
Finally, under the DC-SBM, if the smallest expected degree is growing fast enough, then $P$ converges to $\sP$ in spectral norm (e.g., see \cite{chung2011spectra}). So, estimates of the spectral properties of $\sP$ are similar to the spectral properties of  $P$.  With these facts in mind, we propose estimating the spectral properties of $P$ with the spectral properties of a normalized version of $\hat Q$.
We let $Z \in \{0,1\}^{N \times K}$ be such that $Z_{ij} = 1$ if and only if $z(i)  = j$.
\begin{prop}\label{prop:Q}
If $\sP$ is constructed from the DC-SBM and if $\hat Q$ is computed via a sample from the $(\T,\sP)\textit{-walk on } G$, then
\[\E(\hat Q) =  {\bf B}/m,\]
where $m = 1^T {\bf B} 1$.  
\end{prop}
\begin{prop}  \label{prop:B}
Define  $D_B \in \R^{K \times K}$ to be a diagonal matrix that contains the row sums of ${\bf B}\in \R^{K \times K}$, i.e., $D_B = \diag( {\bf B} \1_K)$, and define $B_L = D_B^{-1/2} {\bf B}D_B^{-1/2}$.    Define $U$ and $\Lambda$ via the eigendecomposition, $B_L = U \Lambda U^T$. 
Define $\beta_\ell^* = \langle y, f_\ell^* \rangle_{\pi^*}$, where $\pi^*$ is the stationary distribution of $\sP$.   Then, (i) the nonzero eigenvalues of $B_L$ are identical to the nonzero eigenvalues of $\sP$, (ii) the columns of 
 \begin{equation} \label{eq:fFormula}
 f^* =\sqrt{m} Z D_B^{-1/2} U
 \end{equation}
 are eigenvectors of $\sP$, and (iii) if $X$ is sampled from $\pi^*$, then 
\begin{equation}\label{eq:betaFormula}
\beta_\ell^* = \E(y(X) f_\ell^*(X)), \mbox{ for } \ell \le K.
\end{equation}
\end{prop}
The proofs of the propositions are given in Section~\ref{sec:technical}. We now introduce our estimator of $\Sigma$ and $\mu$.

\medskip

\noindent
\textbf{SBM-fGLS:}  Using $\tilde z : \T \rightarrow \{1, \dots, K\}$ as an observed partition of the nodes (e.g., by demographic characteristics), the SBM estimator of $\Sigma$  is computed with the following steps.  Each step uses a plug-in estimator using the previously derived formulas.  After the statement of the algorithm, the steps are matched to the motivating equation.

For notational convenience, denote $Y_\tau$, $\tilde z(\tau)$, and $deg(\tau)$ as $y(X_\tau), \tilde z(X_\tau),$ and $deg(X_\tau)$ for each sampled individual $X_\tau$. Moreover, suppose a one-to-one mapping  between the node set of $\T$ and $\{1, \dots, n\}$.

\begin{enumerate}
\item[1)] Compute $\hat Q$ via \eqref{Qdef} using the block memberships $\tilde z(\tau)$.  Define $\hat Q^{(S)} = (\hat Q + \hat Q^T)/2$. This symmetrization ensures the  eigenvalues are real-valued. 
\item[2)] Row and column normalize $\hat Q^{(S)}$, as $\hat Q_L = D_{\hat Q}^{-1/2} \hat Q^{(S)} D_{\hat Q}^{-1/2},$
where $D_{\hat Q} = \diag(\hat Q \1_K) \in \R^{K \times K}$.
\item[3)] Take an eigendecomposition of 
\begin{equation}\label{eq:lamsbm}
\hat Q_L = \hat U \hat \Lambda \hat U^T.
\end{equation}
\item[4)] Compute $\hat f = \hat Z D_{\hat Q}^{-1/2} \hat U$, where $\hat Z \in \{0,1\}^{n\times K}$ contains $\hat Z_{ij} = 1$ iff $\tilde z(i) = j$. 
\item[5)] For $\ell = 1, \dots, K$, compute $\hat \beta_\ell = \frac{1}{n} \sum_\tau Y_\tau \hat f_\ell(\tau)$, where $\hat f_\ell(\tau)$ is the $(\ell, \tau)$ element of $\hat f$.
\item[6)] Compute an estimate of the auto-covariance function as
\[\hat \gamma_{\sbm}(d) = \sum_{\ell = 1}^K  \hat \beta_\ell^2 \hat \Lambda_{\ell \ell}^d.\]
\item[7)] Define $\hat s^2$ to be the sample variance of the $Y_\tau$.  
For $\sigma, \tau \in \T$,
$$\hat \Sigma_{\sigma, \tau}^{\sbm} = 
\left\{\begin{array}{lc}
\hat \gamma_{\sbm}(d(\sigma, \tau)) & \text{if $\sigma \ne \tau$} \\
\hat \gamma_{\sbm}(0) + \hat s^2& \text{if $\sigma  = \tau$,} \end{array}\right.$$
where $\hat s^2$ provides for Tikhonov regularization in $(\hat \Sigma^{\sbm})^{-1}$.
 \item[8)] Define $\hat g \in \R^n$ to solve the system of equations $\hat \Sigma_{sbm} \hat g = \textbf{1}$.
\item[9)] Estimate $\E(Y_\tau)$ with $\sum_{\tau \in \T} \hat g_\tau Y_\tau / \sum_{\tau \in \T} \hat g_\tau$.
\end{enumerate}

Step 1) comes from Proposition \ref{prop:Q}.  Steps 4) and 5) come from  \eqref{eq:fFormula} and \eqref{eq:betaFormula} in Proposition \ref{prop:B}.  Step 6) comes from \eqref{eq:gammadef}.  In all of the plug-in formulas, it is unnecessary to estimate $m$ because we must only specify $\hat \Sigma$ up to a constant of proportionality; this constant appears in both the numerator and denominator of $\hat \mu_{\mathrm{fGLS}}$ in step 9.



\section{Simulations} \label{sec:sim}

\begin{figure*}[t] 
   \centering
\textbf{fGLS estimators  reduce the RMSE when the outcome is hard to estimate}
\newline
  \includegraphics[width=5in]{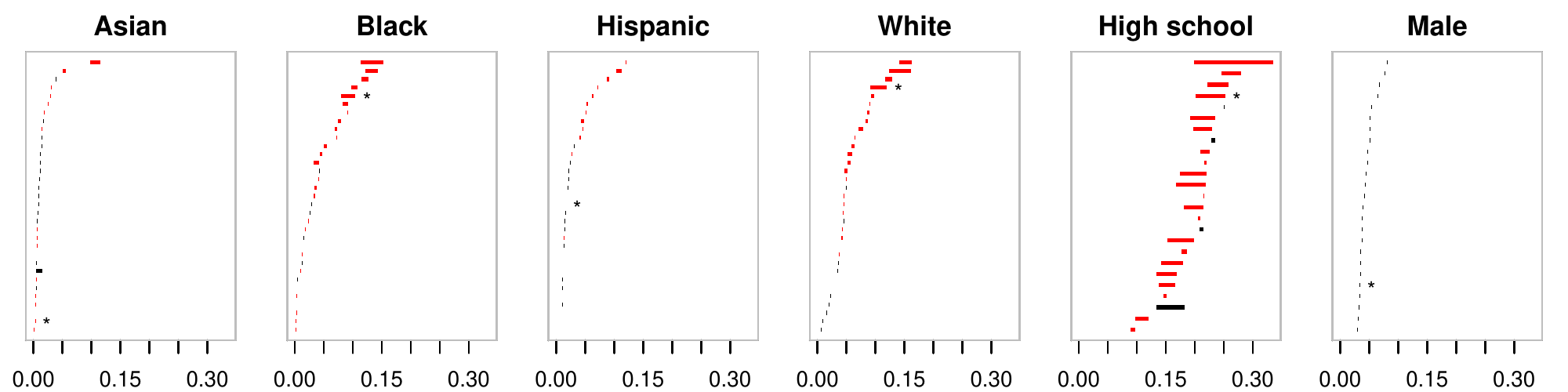} 
   \caption{
These figures present the root mean squared error (RMSE)  for the SBM-fGLS estimator and the Volz-Heckathorn (VH) estimator.  Each panel corresponds to a different outcome $y$.  In each panel, the horizontal axis corresponds to RMSE and the vertical axis corresponds to different schools, ordered by RMSE of the VH estimator. Each line connects the RMSE for the SBM-fGLS to the RMSE of the VH estimator.  If the line is red, then SBM-fGLS has a smaller RMSE.  
  }   \label{fig:rmse}
\end{figure*}

\begin{figure*}[t] 
   \centering
\textbf{These plots diagnose whether fGLS will improve estimation}
\newline
  \includegraphics[width=5in]{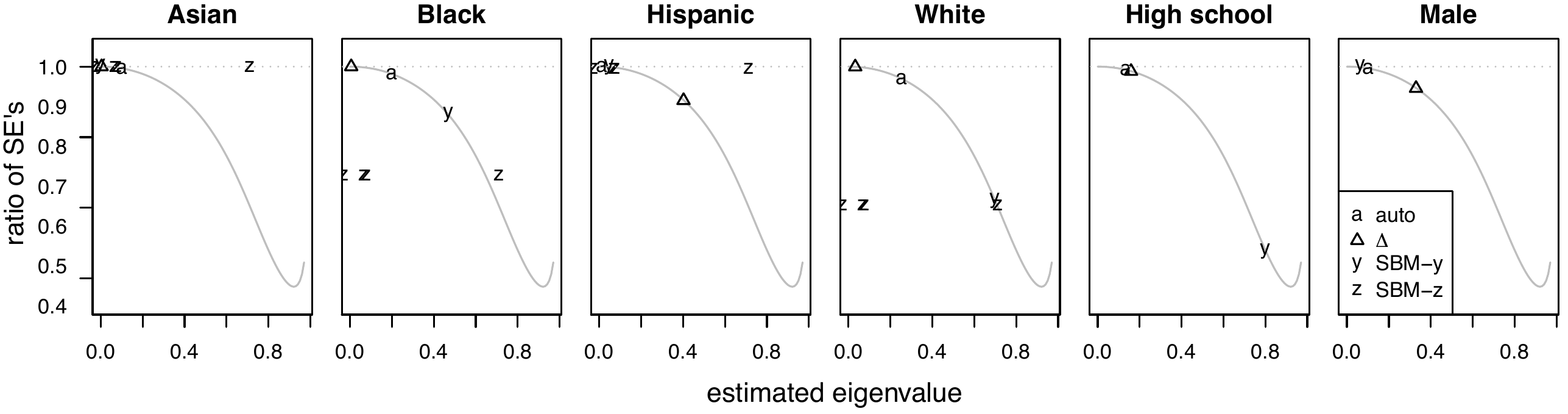} 
   \caption{Each of these diagnostic plots is created from a single sample on the school with the asterisk in Figure \ref{fig:rmse}.  
   We should prefer the fGLS estimators that have a smaller ratio of standard errors (RSE) as defined in the text and displayed on the vertical axis.  
   The \texttt{y} corresponds to the SBM-fGLS estimator that constructs the blocks from the outcome variable of interest.  For the race and ethnicity outcomes, \texttt{z} corresponds to the SBM-fGLS estimator that constructs the blocks with all races and ethnicities observed in the sample.  
   In each plot, there are $(K-1)$-many ${\bf z}$s because SBM-fGLS estimates $K-1$ eigenvalues; each of these $K-1$ points has the same value on the vertical axis. 
   For completeness, this plot includes the rank-two estimators $\hat \mu_{\auto}$ and $\hat \mu_{\Delta}$ that are developed in Section S4 in the SI.  Under the rank-two model, the ratio of SE's is completely determined by the estimated eigenvalue; this is the grey line.  
}   \label{fig:diagnostics}
\end{figure*}

This section compares the SBM-fGLS estimator to the Volz-Heckathorn (VH) estimator via simulation.  Each simulated sample is collected by tracing contacts in social graphs collected in the National Longitudinal Study of Adolescent Health (Add Health).  In the 1994-95 school year, the Add Health study collected a nationally represented sample of adolescents in grades 7-12. The sample covers 84 pairs of middle and high schools in which students nominated of up to five male and five female friends in their middle/high school network \cite{harris2009national}. In this analysis, all contacts are symmetrized to create a social graph.  Then, all graphs are restricted to the largest connected component.    Section~\ref{sec:col} performs a similar simulation on the Colorado Spring Project 90 network \cite{klovdahl1994social}.  These networks were previously studied in \cite{goel2010assessing} and \cite{baraff2016estimating}.  

The RDS process is simulated without replacement.  First, the seed node is selected with probability proportional to node degree.  Then, for each participant $i$, let $R_i$ be an iid random variable with $P(R_i = 0) = 1/6, P(R_i = 1) = 1/3, P(R_i = 2) = 1/3, P(R_i = 3) = 1/6$.  Each participant refers a group of $R_i$ contacts, selected uniformly at random from their contacts whom have not yet participated.  If the participant has fewer than $R_i$  contacts eligible to refer, then the participant refers all of their eligible contacts.  This process continues one referral wave at a time.   The process stops when there are 500 participants.  If the referral process terminates before collecting 500 participants, then the process is restarted.   We restrict the analysis to the 29 different Add Health networks with at least 1000 nodes (i.e. twice the sample size). On each network, we collect 200 different RDS samples.

We construct both SBM-fGLS and VH estimators for the proportion of students of the four largest  race and ethnicity categories (Asian, black, Hispanic, white), the proportion of students in high school, 
and the proportion of students that are male.  The SBM-fGLS estimators for the race and ethnicity categories were constructed with $K=6$ blocks.  
The first four blocks are Asian, black, Hispanic, and white.  Then, the fifth and sixth race/ethnicity categories come from missing and multiple race/ethnicity.  
All other estimators were constructed with $K=2$ blocks, where the blocks were defined by the outcome variable.

In these simulations, the fGLS estimator reweights the outcome $Y$ to adjust for the sampling bias (akin to Remark \ref{remark:unbiasWeighting}), but with an fGLS estimate of the normalizing constant $\E(1/\deg(X))$.  
We use SBM-fGLS to estimate this quantity, using the same blocks as used to construct the estimator for the outcome of interest (e.g. race and ethnicity). See Section~\ref{sec:construction} for a step-by-step construction of the SBM-fGLS estimator.

Figure \ref{fig:rmse} shows the root mean squared error (RMSE) for fGLS and VH estimators. 
Overall, the SBM-fGLS estimator has a substantially smaller RMSE for the hard-to-estimate quantities. 
Each panel in Figure \ref{fig:rmse} has one line with an asterisk.  These lines correspond to the same school, which has both (i) a referral bottleneck between the white and black populations and (ii) a referral bottleneck between the high school and middle school. None of the fGLS estimators model {\em both} bottlenecks,  
yet they perform well. 


\subsection{Diagnostic plot} 
Figure \ref{fig:diagnostics} presents a diagnostic plot to evaluate the fGLS estimators using only data that is observed in a single sample.  This diagnostic plot was created from the first simulated sample taken on the school that has the asterisk in Figure \ref{fig:rmse}.  

The horizontal axis in Figure \ref{fig:diagnostics} gives eigenvalue(s) of $P$ estimated by the fGLS technique.  The vertical axis gives the plug-in estimate for the ratio of standard errors, 
\[\mathrm{RSE}(\hat \Sigma) = \sqrt{
\frac{\widehat{\var}(\hat \mu_{\gls}(\hat \Sigma))} 
{\widehat{\var}(\hat \mu)}
} = 
\sqrt{\frac
{(\1^T \hat \Sigma^{-1} \1)^{-1} }
{n^{-1} \1^T \hat \Sigma \1}}.\]
We should prefer the fGLS estimators that have a smaller ratio.
Estimators with smaller RSE detect more dependencies and thus make further reductions in the variance.  This is discussed in more detail in Section~\ref{sec:diagnostic}.
Notice how the 
fGLS estimators 
have smaller ratios for the outcomes black, white, and high school. These are the outcomes for which fGLS reduces the RMSE in Figure \ref{fig:rmse}.
For Asian, Hispanic, and male, the ratio of SE's is closer to one.  


%


\section{Summary} \label{sec:disc}

This paper derives and studies GLS and fGLS estimators that account for the covariance between samples in a respondent-driven sample.   
Under the Markov model where the covariance between samples is known, Theorems  \ref{thm:main} and \ref{thm:oneEigen} show that the variance of the GLS estimator decays like $O(n^{-1})$.   
To estimate the covariance between samples, we use the fact that the covariance between adjacent samples can be exactly specified in terms of the spectral properties of the Markov transition matrix \cite{verdery2013network, acrds, treevar}.  These essential spectral  properties of the network can be estimated from the observed data under the DC-SBM and the rank-two model.


Section \ref{sec:sim} shows in simulations on the Add Health networks that the fGLS estimates have smaller RMSE than VH estimates for hard to estimate quantities.   This simulation is performed under a more realistic model than the models used in the technical results (Theorems and Propositions \ref{thm:main}, \ref{thm:oneEigen}, \ref{prop:Q}, \ref{prop:B}).  First, the RDS is simulated on social graphs that were recorded in the Add Health study. Second, the sampling is without replacement.  The construction of the SBM-fGLS estimates 
demonstrates that it is not necessary to model the entire network in order to adequately model the correlation between samples. 
In particular, in many of the schools, there are referral bottlenecks on both race/ethnicity and grade level.  It is not necessary to model both bottlenecks, only the one that is related to the outcome of interest. 

The diagnostic plots in Figure \ref{fig:diagnostics} help to determine whether the outcome of interest is correlated in the observed sample.  For quantities that are correlated (e.g. race, ethnicity, and school), Figure \ref{fig:rmse} shows that fGLS estimates drastically reduce the RMSE. 

 
Sections~\ref{sec:simulations} and~\ref{sec:pref} present two additional simulations to investigate the role of (i) sample size, (ii) referral rates, (iii) alignment of the outcome $y$ with the blocks $z$, and (iv) preferential recruitment.  The simulations show that if the outcome of interest correlates or aligns with the underlying structure of the graph and the referral rate is larger than the critical threshold identified in \cite{treevar}, then fGLS estimators can drastically reduce the variability of previous estimators. In some simulations, the fGLS estimators have a smaller RMSE with 500 samples than the VH estimators have with 1000 samples.  While the fGLS estimators are derived under a Markov model, all simulations were performed under a without-replacement (i.e., non-Markovian) model.  In this sense, all of the estimators appear robust to the Markovian model misspecification.  

 
Neither the SBM-fGLS estimator constructed with the outcome of interest, nor $\hat \mu_{\auto}$  require any additional information.  In the simulations, they rarely increase the RMSE and sometimes drastically reduced the RMSE.  Moreover, the diagnostic plots help to identify when these estimators have detected correlation in the sample.  Finally, Section~\ref{sec:pref} shows that these benefits continue to hold under preferential recruitment, where all of the estimators become highly biased.  

\bigskip

\noindent
\textbf{Acknowledgements:} Roch is supported by NSF grants DMS-1149312 (CAREER) and DMS-1614242. Rohe is supported by NSF grant DMS-1612456 and ARO grant W911NF-15-1-0423.

\clearpage

\bibliographystyle{alpha}
\bibliography{TV-references}

\newcommand{\etalchar}[1]{$^{#1}$}
\begin{thebibliography}{PWM{\etalchar{+}}04}

\bibitem[Ait36]{aitken1936iv}
Alexander~C Aitken.
\newblock {IV.}—on least squares and linear combination of observations.
\newblock {\em Proceedings of the Royal Society of Edinburgh}, 55:42--48, 1936.

\bibitem[Ame85]{amemiya1985advanced}
Takeshi Amemiya.
\newblock {\em Advanced econometrics}.
\newblock Harvard university press, 1985.

\bibitem[BMR16]{baraff2016estimating}
Aaron~J Baraff, Tyler~H McCormick, and Adrian~E Raftery.
\newblock Estimating uncertainty in respondent-driven sampling using a tree
  bootstrap method.
\newblock {\em Proceedings of the National Academy of Sciences}, page
  201617258, 2016.

\bibitem[BP94]{benjamini1994markov}
Itai Benjamini and Yuval Peres.
\newblock Markov chains indexed by trees.
\newblock {\em The Annals of Probability}, pages 219--243, 1994.

\bibitem[CR11]{chung2011spectra}
Fan Chung and Mary Radcliffe.
\newblock On the spectra of general random graphs.
\newblock {\em The electronic journal of combinatorics}, 18(P215):1, 2011.

\bibitem[Dur07]{durrett2007random}
Richard Durrett.
\newblock {\em Random graph dynamics}, volume 200.
\newblock Cambridge university press Cambridge, 2007.

\bibitem[Gil11]{gile2011improved}
Krista~J Gile.
\newblock Improved inference for respondent-driven sampling data with
  application to hiv prevalence estimation.
\newblock {\em Journal of the American Statistical Association}, 106(493),
  2011.

\bibitem[GS09]{goel2009respondent}
Sharad Goel and Matthew~J Salganik.
\newblock Respondent-driven sampling as markov chain monte carlo.
\newblock {\em Statistics in medicine}, 28(17):2202--2229, 2009.

\bibitem[GS10]{goel2010assessing}
Sharad Goel and Matthew~J Salganik.
\newblock Assessing respondent-driven sampling.
\newblock {\em Proceedings of the National Academy of Sciences},
  107(15):6743--6747, 2010.

\bibitem[Hec97]{heckathorn1997respondent}
Douglas~D Heckathorn.
\newblock Respondent-driven sampling: a new approach to the study of hidden
  populations.
\newblock {\em Social problems}, 44(2):174--199, 1997.

\bibitem[HHW{\etalchar{+}}09]{harris2009national}
Kathleen~M Harris, Carolyn~T Halpern, Eric Whitsel, Jon Hussey, Joyce Tabor,
  Pamela Entzel, and J~Richard Udry.
\newblock The national longitudinal study of adolescent health: Research
  design.
\newblock {\em Available at h ttp://www. cpc. unc. edu/pr
  ojects/addhealth/design}, 2009.

\bibitem[HLL83]{holland1983stochastic}
P.W. Holland, K.B. Laskey, and S.~Leinhardt.
\newblock Stochastic blockmodels: First steps.
\newblock {\em Social Networks}, 5(2):109--137, 1983.

\bibitem[HT52]{horvitz1952generalization}
Daniel~G Horvitz and Donovan~J Thompson.
\newblock A generalization of sampling without replacement from a finite
  universe.
\newblock {\em Journal of the American statistical Association},
  47(260):663--685, 1952.

\bibitem[KHRR15]{acrds}
Mohammad Khabbazian, Bret Hanlon, Zoe Russek, and Karl Rohe.
\newblock Novel sampling design for respondent-driven sampling.
\newblock {\em In preparation}, 2015.

\bibitem[KN11]{karrer2011stochastic}
B.~Karrer and M.E.J. Newman.
\newblock Stochastic blockmodels and community structure in networks.
\newblock {\em Physical Review E}, 83(1):016107, 2011.

\bibitem[KPW{\etalchar{+}}94]{klovdahl1994social}
Alden~S Klovdahl, John~J Potterat, Donald~E Woodhouse, John~B Muth, Stephen~Q
  Muth, and William~W Darrow.
\newblock Social networks and infectious disease: The colorado springs study.
\newblock {\em Social science \& medicine}, 38(1):79--88, 1994.

\bibitem[LPW09]{levin2009markov}
David~Asher Levin, Yuval Peres, and Elizabeth~Lee Wilmer.
\newblock {\em Markov chains and mixing times}.
\newblock American Mathematical Soc., 2009.

\bibitem[LR15]{li2015central}
Xiao Li and Karl Rohe.
\newblock Central limit theorems for network driven sampling.
\newblock {\em arXiv preprint arXiv:1509.04704}, 2015.

\bibitem[PWM{\etalchar{+}}04]{potterat2004network}
J~Potterat, Donald~E Woodhouse, Stephen~Q Muth, Richard~B Rothenburg, William~W
  Darrow, Alden~S Klovdahl, John~B Muth, et~al.
\newblock Network dynamism: History and lessons of the colorado springs study.
\newblock In {\em Network epidemiology: A handbook for survey design and data
  collection}. Oxford University Press, 2004.

\bibitem[QR13]{qin2013regularized}
Tai Qin and Karl Rohe.
\newblock Regularized spectral clustering under the degree-corrected stochastic
  blockmodel.
\newblock In {\em Advances in Neural Information Processing Systems}, pages
  3120--3128, 2013.

\bibitem[Roh15]{treevar}
Karl Rohe.
\newblock Network driven sampling; a critical threshold for design effects.
\newblock {\em arXiv preprint arXiv:1505.05461}, 2015.

\bibitem[RWP{\etalchar{+}}95]{rothenberg1995social}
Richard~B Rothenberg, Donald~E Woodhouse, John~J Potterat, Stephen~Q Muth,
  William~W Darrow, and Alden~S Klovdahl.
\newblock Social networks in disease transmission: the colorado springs study.
\newblock {\em NIDA research monograph}, 151:3--19, 1995.

\bibitem[SH04]{salganik2004sampling}
Matthew~J Salganik and Douglas~D Heckathorn.
\newblock Sampling and estimation in hidden populations using respondent-driven
  sampling.
\newblock {\em Sociological methodology}, 34(1):193--240, 2004.

\bibitem[VH08]{volz2008probability}
Erik Volz and Douglas~D Heckathorn.
\newblock Probability based estimation theory for respondent driven sampling.
\newblock {\em Journal of Official Statistics}, 24(1):79, 2008.

\bibitem[VMBM13]{verdery2013network}
Ashton~M Verdery, Ted Mouw, Shawn Bauldry, and Peter~J Mucha.
\newblock Network structure and biased variance estimation in respondent driven
  sampling.
\newblock {\em arXiv preprint arXiv:1309.5109}, 2013.

\bibitem[WHO13]{WHO}
{\em Introduction To HIV/AIDS And Sexually Transmitted Infection Surveillance
  Module 4: Introduction to Respondent-drive Sampling.}
\newblock World Health Organization \& UNAIDS, 2013.
\newblock http://applications.emro.who.int/dsaf/EMRPUB\_2013\_EN\_1539.pdf.

\bibitem[WHS{\etalchar{+}}15]{white2015strengthening}
Richard~G White, Avi~J Hakim, Matthew~J Salganik, Michael~W Spiller, Lisa~G
  Johnston, Ligia Kerr, Carl Kendall, Amy Drake, David Wilson, Kate Orroth,
  et~al.
\newblock Strengthening the reporting of observational studies in epidemiology
  for respondent-driven sampling studies:``strobe-rds'' statement.
\newblock {\em Journal of clinical epidemiology}, 68(12):1463--1471, 2015.

\bibitem[WRP{\etalchar{+}}94]{woodhouse1994mapping}
Donald~E Woodhouse, Richard~B Rothenberg, John~J Potterat, William~W Darrow,
  Stephen~Q Muth, Alden~S Klovdahl, Helen~P Zimmerman, Helen~L Rogers, Tammy~S
  Maldonado, John~B Muth, et~al.
\newblock Mapping a social network of heterosexuals at high risk for hiv
  infection.
\newblock {\em Aids}, 8(9):1331--1336, 1994.

\end{thebibliography}

\appendix

\section{Simulation on network from Colorado Springs\\ Pro\-ject 90}\label{sec:col}

Between 1988 and 1992, a study funded by the United States Centers for Disease Control and Prevention mapped the HIV transmission network among injection drug users and sex workers in Colorado Springs.  We use the largest connected component of this graph which contains 4,430 individuals and 18,407 edges \cite{klovdahl1994social,woodhouse1994mapping, rothenberg1995social,potterat2004network}.  We simulate the RDS process the same as was done on the Add Health networks in the main text, with one exception.  Because this network is much larger, we also include a simulation for a sample size of 1000.  We use the SBM-fGLS estimator where the blocks are constructed with each outcome variable.  For each trait that we estimate, there are $K=2$ blocks.

\begin{figure}[h] 
	\centering
	\includegraphics[width=5in]{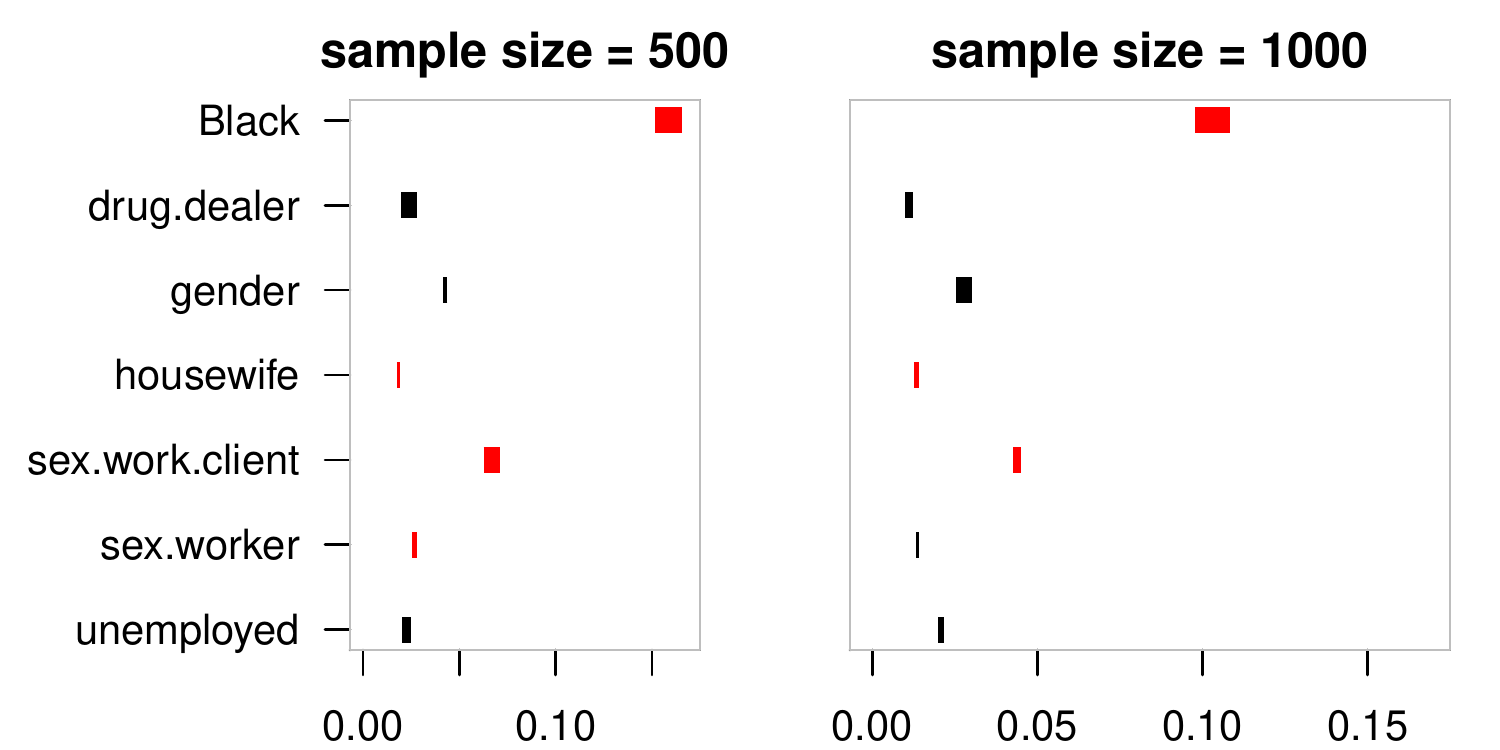} 
	\caption{These figures presents the root mean squared error (RMSE) for the SBM-fGLS estimator and the Volz-Heckathorn (VH) estimator on seven different traits and for two different sample sizes. The panel corresponds to a sample size of 500 and the right panel corresponds to a sample size of 1000.  In each panel, the horizontal axis corresponds to RMSE and the vertical axis corresponds to different traits, ordered alphabetically.  Each line connects the RMSE for the SBM-fGLS to the RMSE of the VH estimator. If the line is red, then SBM-fGLS has a smaller RMSE.
	}   \label{fig:colRMSE}
\end{figure}

There are several traits recorded on each individual.  If a value was missing for an individual, then we presume that the individual does not have that trait.  We include all traits which  for which at least 5\% of the population has that trait.  Previous research in \cite{goel2010assessing} has show that it is particularly difficult to estimate the proportion of individuals that are black.  This is again confirmed in Figure \ref{fig:colRMSE}.  Moreover, Figure \ref{fig:colRMSE} shows that fGLS reduces the RMSE for this quantity.  

\begin{figure}[h] 
	\centering
	\includegraphics[height=6in]{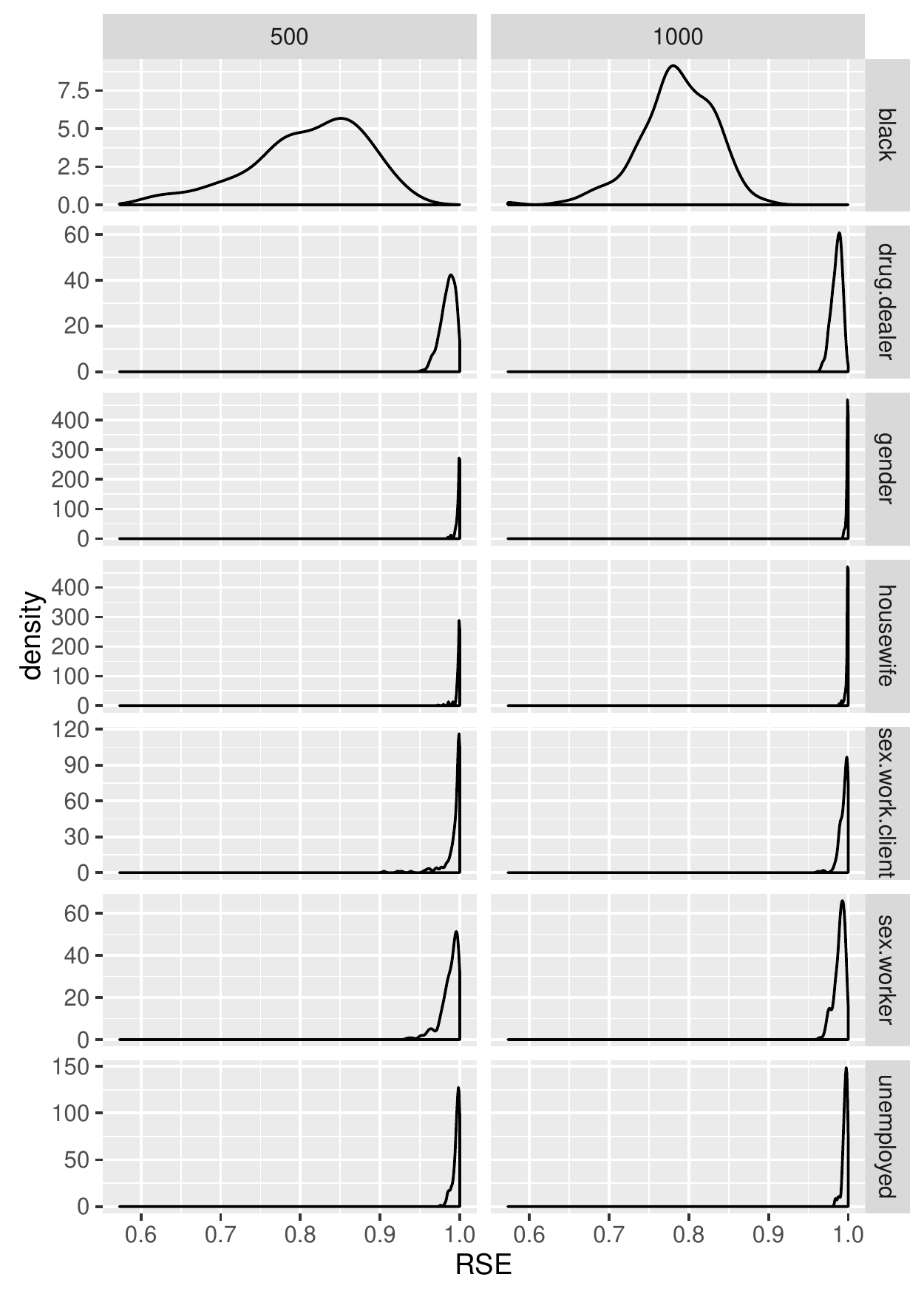} 
	\caption{These figures give the distribution of the $\rse$ over the 200 simulated samples.  The left column corresponds to the sample size of 500 and the right panel corresponds the sample size of 1000.  Each row corresponds to a different trait.  All estimators are SBM-fGLS constructed from the outcome of interest.  $\rse$ overstate how much fGLS reduces the standard error,  as discussed in Section \ref{sec:diagnostic}.
	}  \label{fig:colDiagnostics}
\end{figure}

Figure \ref{fig:colDiagnostics} presents a histogram of the $\rse$ over the 200 simulated samples.  $\rse$ was defined in the main text (and is also in Equation \eqref{eq:rse} below).  This is the value along the vertical axis in the diagnostic plot.  Figure \ref{fig:colDiagnostics} shows that $\rse$ is much smaller when estimating the proportion of individuals that are black.   For the other traits, the RSE is much closer to one.  As such, the diagnostic plot again selects the variable for which it is most beneficial to use fGLS.  

Figure \ref{fig:colDiagnostics} shows that the actual values of $\rse$ overstate how much fGLS reduces the standard error.  As discussed in Section \ref{sec:diagnostic},  $\rse$ should not be considered as an estimator because it is biased.  However, as a diagnostic quantity, it successfully identifies the variables that benefit from fGLS.

\section{Simulation with simulated network}
\label{sec:simulations}
\begin{figure}[htbp] 
	\centering
	\includegraphics[width=3in]{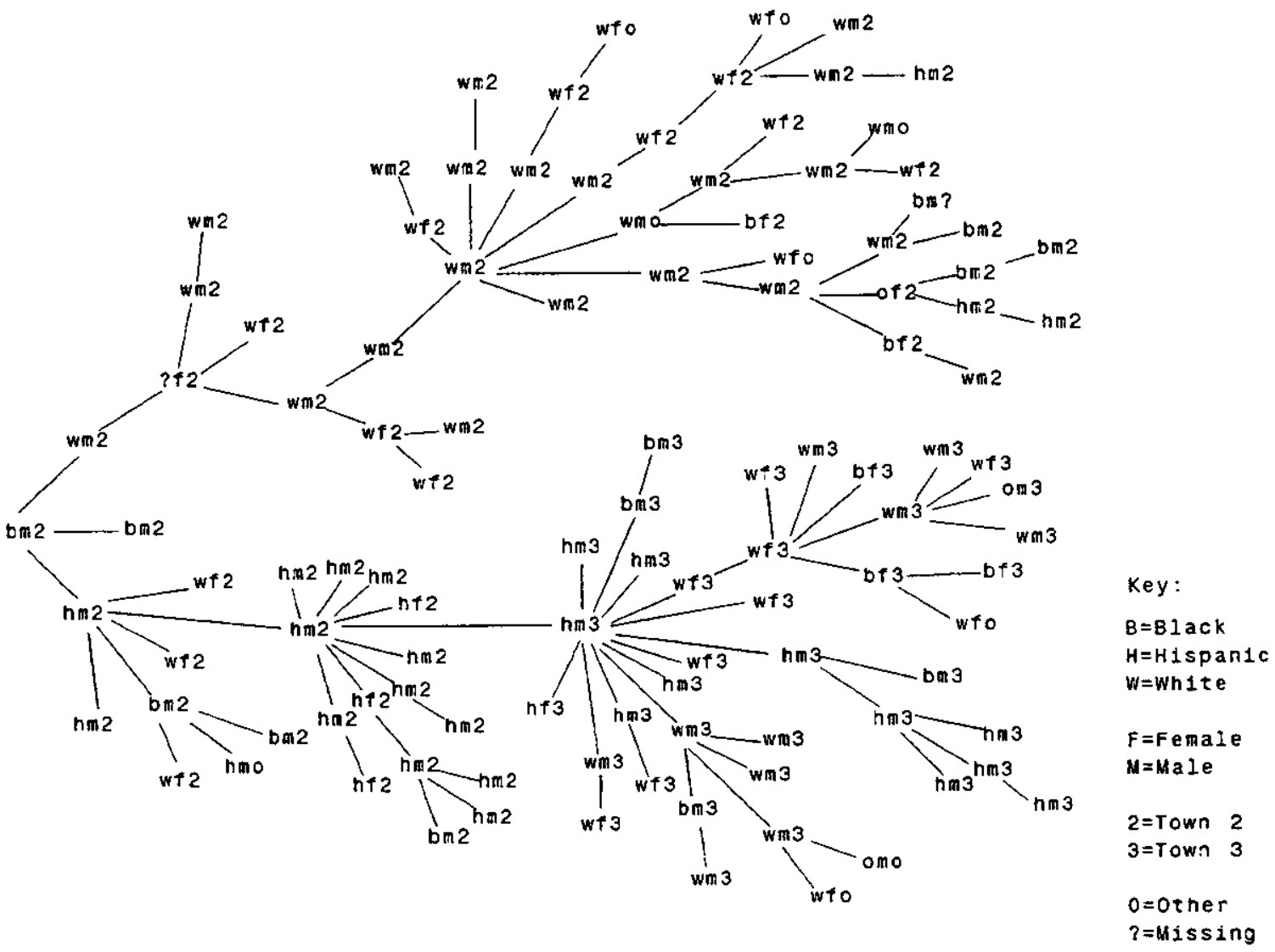} 
	\caption{Reprinted referral tree from \cite{heckathorn1997respondent}.} \label{fig:heck}
\end{figure}

Figure \ref{fig:heck} is a reprint of Figure 1 from \cite{heckathorn1997respondent}.  It shows the referral tree which reached 112 injection drug users in Eastern Connecticut of various races, genders, and towns.  In addition to reporting $\T$, the figure also reports the race, gender, and town of each of these 112 participants.  
In order to create a more realistic simulation, but also investigate and control key elements of the network, we use the data in this figure to pick 
(i) parameter settings to generate the underlying social network $G$ and (ii)
parameter settings to generate a referral tree $\T$.

The tree in  Figure \ref{fig:heck} gives three demographic measurements on each individual.  Hence,  we constructed three different versions of the matrix $\hat Q$, one for each of the three demographic variables.  Among these three matrices, race produced the most divisive bottleneck, as measured by the second largest eigenvalue of the matrix $\hat Q_L$ from step (2) of the SBM-motivated algorithm above. 
Table \ref{tab:Q} gives the matrix $n \hat Q$ for the partition on race.   
The simulations below use this $\hat Q$ to parameterize the simulation of $G$ (in view of Proposition~\ref{prop:Q}).

\begin{table}
	\centering
	\caption{Matrix $n \hat Q$ counts transitions between {\em B, W, H}.} \label{tab:Q}
	
	\begin{tabular}{l|ccc}
		& \em B&\em W&\em H\\
		\midrule
		\em B &5&5&2\\
		\em W &7&46&1\\
		\em H&4&8&28\\
	\end{tabular}
\end{table}
\vspace{.2in}

Everything that follows uses the symmetrized version of $\hat Q$, $(\hat Q + \hat Q^T)/2$.
To make the expected degree equal to $30$, set ${\bf B} =30 N \hat Q$.
The proportion of nodes in each of the three blocks was chosen to be proportional to the row sums of $\hat Q_L$ (as defined in Step 2 of the SBM estimation of $\Sigma$).  
The resulting block proportions $(p_b, p_w, p_h) \approx (.13, .33, .53)$ are close to the empirical proportions $(.17, .29, .54)$.    
We set $\theta_i = .3 + w_i,$ where $w_i $ are independent $\mathrm{Gamma}(\mathrm{shape} = 200, \mathrm{rate} = 300)$.  The .3 prevents isolated nodes (i.e., nodes of degree zero), while the Gamma distribution ensures some degree heterogeneity within blocks.  After this sampling, the $\theta_i$s in the same block were scaled to sum to one (consistently with (13)). 
Using these settings, an underlying network $G$ was sampled from a DC-SBM with $N = 20,000$ nodes.


We examined three different feature vectors, each with $\mu_{\true} \approx .66$: 
$y^{(a)}$ is aligned with $z$, 
$y^{(c)}$ is correlated with $z$, and 
$y^{(u)}$  is uncorrelated with $z$.  Specifically, for each node $i \in G$,
\[ 
y^{(a)}(i) = 
\left\{\begin{array}{ll}1 & z(i) = \mathrm{B} \\1 & z(i) =\mathrm{W} \\0 & z(i) = \mathrm{H}\end{array}\right., 
\quad
y^{(c)}(i) \sim
\left\{\begin{array}{ll}$Ber(.7)$ & z(i) =\mathrm{B} \\$Ber(.1)$ & z(i) = \mathrm{W} \\$Ber(.9)$ & z(i) = \mathrm{H}\end{array}\right.
\]
and $y^{(u)}(i) \sim \mbox{Ber(.66)}$, where Ber$(p)$ is the Bernoulli distribution with probability  $p$.

The referral tree in \cite{heckathorn1997respondent} contains only  112 participants.  
To investigate other sample sizes, it was necessary to generate larger referral trees (i.e., ones with more samples).   We modified the empirical offspring distribution from the tree on 112 participants and simulated Galton-Watson trees, where each individual refers an i.i.d.~number of participants.  
There are two simulation settings for two different referral trees.  The two trees are Galton-Watson trees generated from two different offspring distributions.  The first offspring distribution is created by removing all of the zeros from the empirical offspring distribution.  It has expected value $\approx 2.36$.  
The second offspring distribution is created by reintroducing 15 zeros.  It has expected value $\approx 1.78$.  
The Markov matrix $P$ corresponding to the underlying network $G$ above has  $\lambda_2 \approx .73$. Hence, when the expected value of the offspring distribution exceeds $1/\lambda_2^2 \approx 1.88$, the variance of previous estimators does not decay like $n^{-1}$ under the \tp \cite{treevar}.  As such, one of the simulated trees exceeds this threshold and one does not.  Importantly, the simulation below is not a \tpns.  Instead, it samples without replacement, as described in the next paragraph.

The simulations were performed by simulating one network of size $N = 20,000$ and two referral trees of size $n=1000$ (one tree for each of the two different offspring distributions). Then, with each of the two trees, we simulated 300 different RDS samples as follows. The seed node was selected uniformly.  Then, each participant referred a group of friends uniformly at random from their friends that had not yet participated (Section S1 performs the simulation with preferential recruitment). This sample is not a \tp because it samples \textit{without replacement}.
If a participant did not have enough friends to make their referrals, then the process was restarted. 
This happened in three of 603 total samples.  

We evaluate the performance of the estimators for three different values of the sample size, $n = 100, 500, 1000$.   To create samples of $n=100$ and $n=500$, we used the first $100$ and $500$ samples in the trees of $n=1000$.

We compared six different estimators.  The first two estimators  do not adjust for the covariance: the Volz-Heckathorn estimator $\hat \mu_{\vh}$ and the sequential sampling estimator $\hat \mu_{\ss}$ \cite{gile2011improved}. Computing $\hat \mu_{\ss}$ requires an estimate of the population size which is set to be the true value $N=20,000$.   
The other four estimators are fGLS estimators.  
The fGLS estimators $\hat \mu_{\auto}$ and $\hat \mu_{\Delta}$ are based upon the rank-two model, which is misspecified with the outcome $y^{(c)}$.  These estimators are described in the SI (Section SI2).  
The first SBM-fGLS estimator $\hat \mu_{Y}$ is based upon a misspecified Stochastic Blockmodel for the underlying graph; it estimates the blocks with the  outcome  $y^{(a)}, y^{(c)},$ or $y^{(u)}$ (i.e. HIV status), depending on the simulation setting. Said another way, $\tilde z = Y$. 

%
%
The second SBM-fGLS estimator $\hat \mu_{Z}$ uses the correctly specified Stochastic Blockmodel.  It sets $\tilde z = z$ (i.e., the blocks were correctly specified with race).   
All of these fGLS estimators use the Volz-Heckathorn weights, $\hat \pi$.
The fGLS estimators $\hat \mu_{\auto}, \hat \mu_{\Delta},$ and $\hat \mu_{Y}$ are based upon misspecified models for $\Sigma$.  Both of the SBM-fGLS estimators use the Volz-Heckathorn sampling weights, where the normalizing constant $\E(1/deg(X))$  is estimated via fGLS.    This matches the construction of the estimators in the main paper and is more fully described in Section \ref{sec:construction}.

\begin{figure}[h] 
	\centering
	\includegraphics[width=5in]{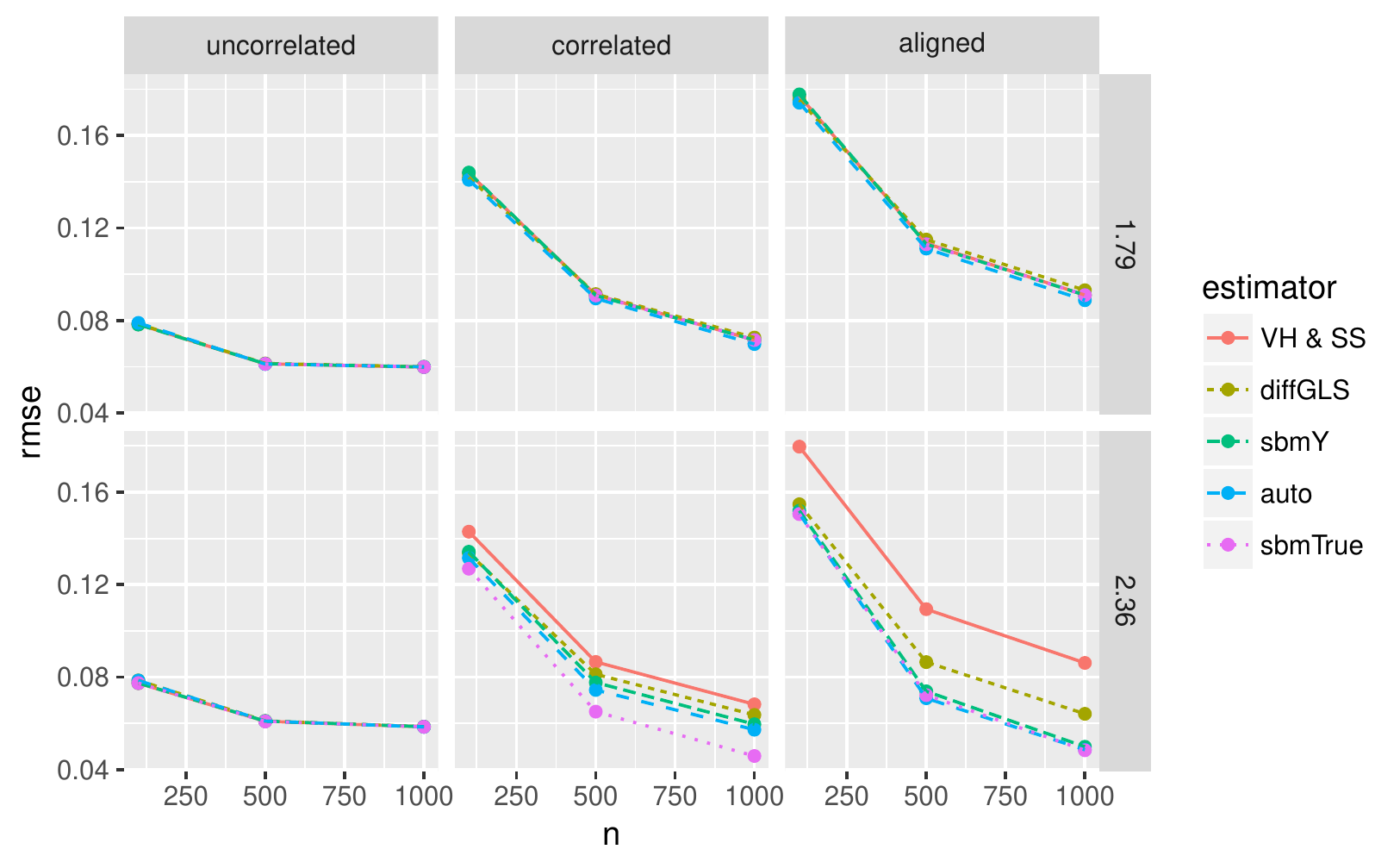} 
	\caption{The three columns correspond to increasing levels of homophily; in the left column HIV status is uncorrelated with the block memberships ($y^{(u)}$) block and
		in the right column nodes HIV is aligned with block memberships ($y^{(a)}$).  The first row corresponds to a slower referral rate.  The second row corresponds to a faster referral rate.  In each panel, the x-axis gives three different sample sizes and the y-axis gives the root mean squared error (RMSE).  Each line corresponds to a different estimator.  In the bottom right panel, the fGLS estimators with 500 samples have smaller RMSE than the standard estimators with 1000 samples. With 1000 samples, the RMSEs of the fGLS estimators (excluding $\hat \mu_{\Delta}$) are nearly half the RMSEs of $\hat \mu_{VH}$ and $\hat \mu_{SS}$. In this figure, SS and VH are indistinguishable.}   \label{fig:rmse-si}
\end{figure}

%
%

Figure \ref{fig:rmse} shows that in certain settings, the fGLS estimators have a smaller the root mean squared error (RMSE) than the VH and SS estimators.  In the top row of plots (corresponding to the smaller referral rate) the fGLS estimators do not improve upon the previous estimators because the simulation is below the critical threshold identified in \cite{treevar}. As such, $\var(\hat \mu_{VH})$ should converge at rate $O(1/n)$.  
Importantly, the fGLS estimators do not increase the RMSE in this regime.  
In the bottom right panel fGLS drastically reduces the RMSE of the previous estimators.  This
corresponds to the simulation setting where the outcome is aligned with the blocks (i.e. $y^{(a)}$) and the  largest referral rate.  
In this panel, all the fGLS estimators with 500 samples have a smaller RMSE than the VH and SS estimators with 1000 samples. Across all other settings, the fGLS estimators never have a larger RMSE than the VH and SS estimators; importantly, this holds even when the fGLS estimators are based upon misspecified models.  In this simulation, we are not sampling a very large proportion of the population.  So, the results for $\hat \mu_{\ss}$ and $\hat \mu_{\vh}$ are indistinguishable in Figure \ref{fig:rmse-si}.  


\section{The benefits of fGLS are robust to preferential recruitment} \label{sec:pref}

In the simulations from the last section, the correlation between samples comes from the fact that friends in the underlying population are likely to have the same value of $y(i)$.  This is called ``homophily.''  Another potential source of correlation that the previous simulations do not include is preferential recruitment, where participants are more likely to refer a friend that shares the same value of $y(i)$.  Recall that in the definition of $P$ in Eq. (2) in the main text,
the probability of a referral from $i$ to $j$ is proportional to $w_{ij}$. Preferential recruitment can be modeled in the \tp by allowing for edges in the graph to have heterogeneous $w_{ij}$s; refer to such a \tp as a \textit{preferential} \tpns.

In the previous simulations, the network was unweighted, $w_{ij}=1$ for all $(i,j)\in E$. To examine the effects of preferential recruitment, we performed another simulation using the same network, with one modification: if $(i,j) \in E$  with $z(i) = z(j) = h$, then $w_{ij} = 10$.  With this, samples in  block $h$ were much more likely to refer other samples in block $h$.  This simulation used the same $\T$s from the simulation in Figure~\ref{fig:diagnostics} in the main text.  
Similarly, this simulation was performed without replacement; the weights $w_{ij}$ were used to sample participants $i$'s friends whom have not yet participated.  Only two of 602 simulations were discarded because a participant did not have enough friends to refer.

\begin{figure}[h] 
	\centering
	\includegraphics[width=5in]{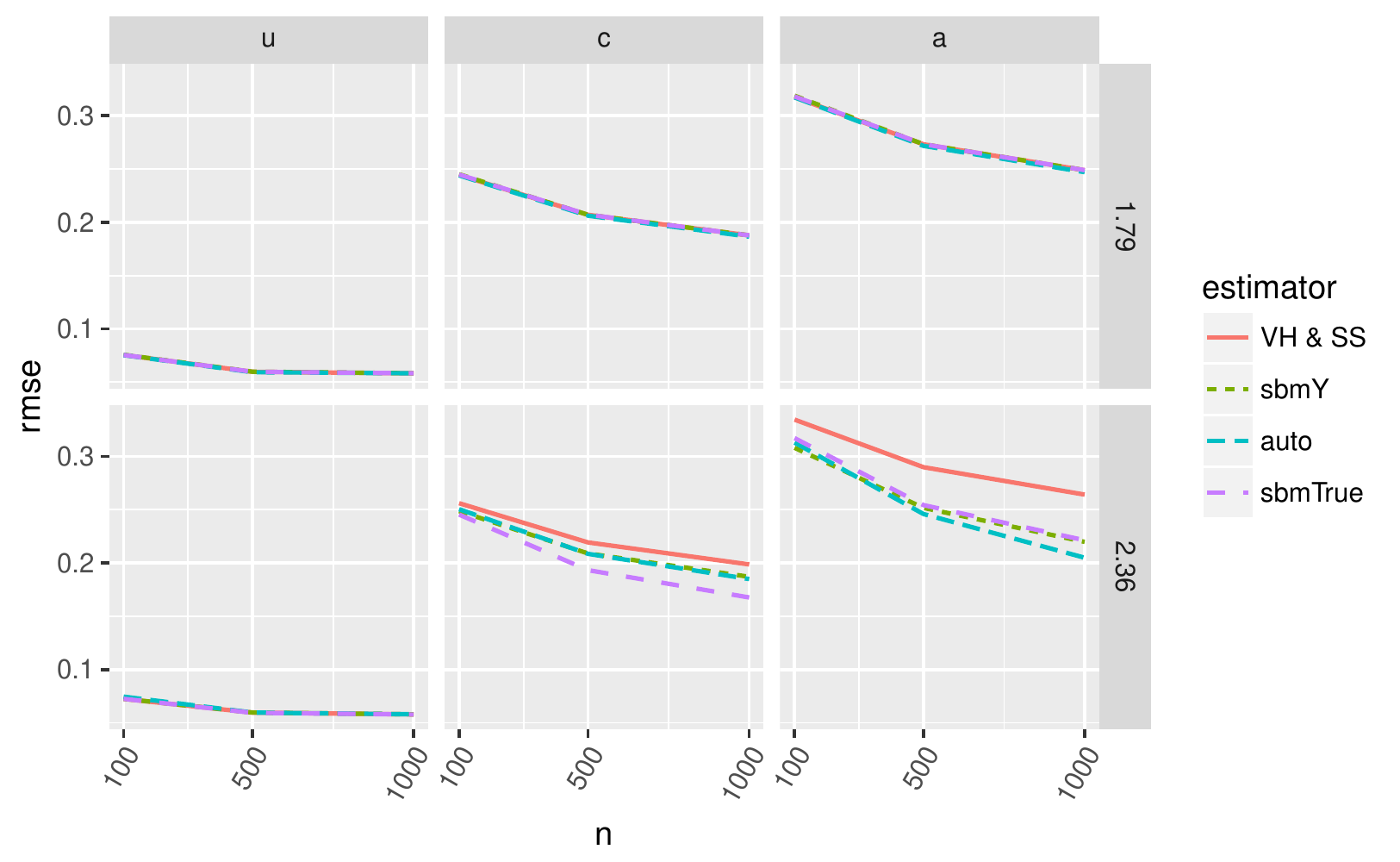} 
	\caption{This figure repeats the simulation from Figure~\ref{fig:diagnostics}, except with the preferential recruitment described in Section \ref{sec:pref}. 
		When the features are correlated or aligned with $z$, then preferential recruitment  makes all of the estimators biased, resulting in an increased RMSE compared to Figure~\ref{fig:diagnostics}.  The fGLS estimators still have a smaller RMSE than $\hat \mu_{VH}$ and $\hat \mu_{SS}$ because they continue to have a smaller variance.   }
	\label{fig:pref}
\end{figure}

Under the \textit{preferential} \tpns, all of the  estimators considered in this paper are biased because they used ``the number of friends'' instead of  $\deg(i) = \sum_j w_{ij}$.   As such, Figure \ref{fig:pref} shows that
preferential recruitment drastically increases the RMSE of all estimators
when the outcome is aligned or correlated with $z$ (i.e., the source of the preference).  



While preferential recruitment drastically increases the RMSE, the fGLS estimators still have a smaller RMSE than $\hat \mu_{\vh}$ and $\hat \mu_{\ss}$.  This highlights how the benefits of fGLS are robust to preferential recruitment.   The reason is that when there is preferential recruitment within block $h$, there is more correlation between samples of $y^{(a)}$ and $y^{(c)}$.  fGLS adjusts for correlation that results from both  homophily and preferential recruitment.




\section{fGLS for the rank-two model} \label{sec:rank2est}

\subsection{\textit{auto} estimator}

Under the rank-two assumption, the auto-covariance at lag $t$ is $\gamma(t) = \beta^2 \lambda^t$.
So,
\[\gamma(0) = \beta^2 \quad \mbox{and} \quad \gamma(1) = \beta^2 \lambda.\]
These equations motivate the estimator $\hat \mu_{\auto}$. 

Let the set $D_k$ contain the node pairs that are distance $k$ apart in $\T$,
\[D_k = \{(\sigma, \tau) : d(\sigma, \tau) = k\}.\]
Define the sample auto-covariance as 
\begin{equation}\label{eq:np}
\hat \gamma_{\np, {\bf m}}(k) = \frac{1}{|D_k|} \sum_{(\sigma, \tau) \in D_k} (Y_\sigma - {\bf m})(Y_\tau - {\bf m}),
\end{equation}
where ${\bf m}$ is some initial estimate of $\mu$.  

%

Denote $\hat \mu_{\auto}({\bf m})$ as the fGLS estimator that comes from the following steps. 
\begin{enumerate}
	\item Use $\hat \gamma_{\np, {\bf m}}$ to compute 
	\begin{equation}\label{eq:lamauto}
	\hat \beta^2 = \hat \gamma_{\np, {\bf m}}(0)    \quad \mbox{and} \quad \hat \lambda = \frac{\hat \gamma_{\np, {\bf m}}(1)}{\hat \gamma_{\np, {\bf m}}(0)}.
	\end{equation}
	\item Compute $\hat \gamma_{\auto}(t) = \hat \beta^2 \hat \lambda^t$.
	\item Define $\hat \Sigma^{({\bf m})} \in \R^{n \times n}$ such that $\hat \Sigma^{({\bf m})}_{\sigma, \tau} = \hat \gamma_{\auto}(d(\sigma, \tau))$, where $d(\sigma, \tau)$ is the graph distance between $\sigma$ and $\tau$ in $\T$.
	\item $\hat \mu_{\auto}({\bf m}) = \hat \mu_{\gls}(\hat \Sigma^{({\bf m})})$, as defined in Eq. (12) in the main text.
\end{enumerate}

The estimator $\hat \mu_{\auto}$ is defined as
\begin{equation} \label{eq:muauto}
\hat \mu_{\auto} = {\arg\min}_{{\bf m} \in [0,1]} \left|\hat \mu_{\auto}({\bf m}) - {\bf m}\right|.
\end{equation}
This is a generalization of the fixed point, ${\bf m}^* = \hat \mu_{\auto}({\bf m}^*)$; such a fixed point does not always exist in our simulations. The range $[0,1]$ was used in the simulations of this paper because $y_i \in \{0,1\}$ for all nodes $i$. 
If $Y_\tau = Y_0$ for all $\tau \in \T$, define $\hat \mu_{\auto} = Y_0$. 


\subsection{$\hat \mu_{\Delta}$}
The estimator $\hat \mu_{\auto}$ is complicated by the fact that  it requires searching over values of $\bf m$.  An alternative approach does not require this search.
Define the function $\Delta: \{0, 1, 2, \dots\} \rightarrow \R$ such that $\Delta(d(\sigma, \tau)) = \E(Y_\sigma - Y_\tau)^2$.  Under the rank-two assumption, if $d(\sigma, \tau) =k$, then
\[\Delta(k) = 2\E(Y_\sigma - \mu)^2  - 2 \E(Y_\sigma - \mu)  (Y_\tau - \mu) = 2\beta^2 -2\,\cov(Y_\sigma,Y_\tau) = 2\beta^2( 1- \lambda^k). \]
So, 
\[\frac{\Delta(2)} {\Delta(1)}  = \frac{2\beta^2( 1- \lambda^2)}{2\beta^2( 1- \lambda)} = \frac{( 1- \lambda)( 1+ \lambda)}{ 1- \lambda} = 1+\lambda.\]
Estimating $\Delta(1)$ and $\Delta(2)$ does not require an estimate of $\mu$, namely,
\begin{equation}\label{eq:DeltaHat}
\hat \Delta(k) = \frac{1}{|D_k|} \sum_{(\sigma, \tau) \in D_k} (Y_\sigma -Y_\tau)^2, \quad\quad D_k = \{(\sigma, \tau) : d(\sigma, \tau) = k\}.
\end{equation}
This yields a plug-in estimator 
\[\hat \lambda_{\Delta: \mbox{plug-in}} = \frac{\hat \Delta(2)} {\hat \Delta(1)}-1,\]
which often resulted in estimates of $\mu$ that were far outside the range of values in $y$ (simulations not shown).
Rearranging terms and then adding $n^{-1/2}$ to the denominator for Laplace smoothing yields the estimator that we use,
\begin{equation}\label{eq:lamdelta}
\hat \lambda_\Delta = 
\frac{\hat \Delta(2) - \hat \Delta(1)} 
{\hat \Delta(1)+n^{-1/2}}.
\end{equation}
Then, 
\[\hat \Sigma_{\sigma, \tau}^{\Delta}  = \hat \lambda_\Delta^{d(\sigma, \tau)}.\]
It is not necessary to estimate $\beta^2$ because in the computation of 
$$\hat \mu_{\mathrm{fGLS}} = (\1 \hat \Sigma^{-1} \1)^{-1} \1 \hat \Sigma^{-1} Y,$$ any estimate for $\beta^2$ cancels out; it is only necessary to estimate $\hat \Sigma_{\sigma, \tau}^{\Delta}$ up to a proportionality constant.  Denote the resulting estimator 
\[\hat \mu_{\Delta} = \hat \mu_{\gls}( \hat \Sigma^\Delta).\]

\noindent
\textbf{Using $\hat \gamma_{\np, {\bf m}}$ for a fully non-parametric estimate of $\gamma$:
}
Because $\T$ is growing geometrically,
we must only estimate $O(\log n)$ terms (with $n$ samples) in the covariance function $\gamma$ to estimate $\Sigma$.
We were unable to get this approach to perform well in simulations (not shown). 

\section{Technical results}
\label{sec:technical}

\subsection{Proof of Theorem~\ref{thm:oneEigen}}

\begin{proof} From Eq. (4) and Eq. (5) in the main text, $\var(\hat \mu_{\gls}) = (\1^T \Sigma^{-1} \1)^{-1}$.
	Under the rank-two assumption, 
	the inverse of $\Sigma$ has a closed form solution that is sparse. 
	Let $\sigma \sim \tau$ denote an edge in $\T$ between $\sigma$ 
	and $\tau$ and let $\deg(\sigma)$ be the degree of $\sigma$ in $\T$ (ignoring edge direction).  
	It can be seen that the inverse of the covariance matrix $[\Sigma_{\sigma,\tau}]_{\sigma,\tau} = [\gamma(d(\sigma, \tau))]_{\sigma,\tau}$ with
	$\gamma(d) =  \beta^2 \lambda^{d}$
	satisfies
	\begin{equation}\label{eq:sigInverse}
	\beta^{2}(1-\lambda^2) \Sigma^{-1}_{\sigma, \tau} = \left\{\begin{array}{cc}
	1 + \lambda^2(\deg(\sigma) -1) & \sigma = \tau \\
	-\lambda & \sigma \sim \tau \\
	0 & \text{o.w.}\end{array}\right.
	\end{equation}
	
	To compute $\1^T \Sigma^{-1} \1$, first compute $\Sigma^{-1} \1$,
	\begin{eqnarray*}
		\beta^{2}(1-\lambda^2) [\Sigma^{-1} \1]_\sigma 
		&=& 1+ \lambda^2(\deg(\sigma) -1) - \deg(\sigma)\lambda \\
		&=& (1-\lambda)^2 - \deg(\sigma)\lambda(1-\lambda). 
	\end{eqnarray*}
	Factoring out $(1-\lambda)$ on the left- and right-hand sides, then rearranging terms yields
	\begin{equation}\label{eq:solveLinearSystem}
	\beta^{2}(1+ \lambda) [\Sigma^{-1} \1]_\sigma 
	= 1 - \lambda(\deg(\sigma)-1).
	\end{equation}
	This expression provides a closed-form solution to the linear system $\Sigma x = \1$, 
	which is essential to computing the GLS estimator.
	For an undirected and connected tree, note that $\sum_{\sigma\in \T} \deg(\sigma) = 2(n-1)$.
	So, summing over all nodes $\sigma$, we get
	\begin{eqnarray} 
	\beta^{2} (1+ \lambda) \1^T \Sigma^{-1} \1 
	&=& \sum_{\sigma\in \T} \left[1 - \lambda(\deg(\sigma)-1) \right]\\
	&=& n  - \lambda \left( 2(n-1)-n\right) \\
	&=& n \left(1 - \lambda \left( 1 -\frac{2}{n}\right)\right). \label{eq:1sig1}
	\end{eqnarray}
	Rearranging terms gives the result.
\end{proof}

\subsection{Proof of Theorem~\ref{thm:main}} \label{sec:mainproof}

\begin{proof}
	Because the GLS estimator minimizes the variance over all linear estimators, it suffices to construct an explicit
	linear estimator with asymptotic variance $O(1/n)$.
	
	Let $\mathcal{K} = \{\ell\,:\, \beta_\ell \neq 0\}$
	and $K = |\mathcal{K}|$.
	Let $\T$ have $H$ levels, where $H = \Theta(\log n)$
	by the assumption that $\T$ is a complete binary tree.
	Let $\tilde{n}$ be the number of nodes on level
	$H-K$ starting from the root and note that 
	we have $\tilde{n} = \Theta(n)$ as $n \to +\infty$. Denote
	by $\tau_{1,b}$, $b = 1,\ldots,\tilde{n}$,
	the nodes on level $H-K$. For each $b$,
	choose one child $\tau_{2,b}$ of $\tau_{1,b}$,
	one child $\tau_{3,b}$ of $\tau_{2,b}$,
	and so on until the leaves are reached.
	To simplify the notation, we let
	$Z_{a,b} = X_{\tau_{a,b}}$ and 
	$Y_{a,b} = y(Z_{a,b})$ for $a = 1,\ldots,K$ 
	and $b = 1,\ldots,\tilde{n}$.
	
	We consider a linear estimator of the form
	$$
	\Gamma 
	= \frac{1}{\tilde{n}} \sum_{b=1}^{\tilde{n}} \Gamma_b,
	$$
	where
	$$
	\Gamma_b
	= \sum_{a=1}^{K} \gamma_a Y_{a,b},
	$$
	for constants $\gamma_a$, $a=1,\ldots,K$,
	to be determined. 
	The $\gamma_a$s are chosen so that
	$\sum_a \gamma_a = 1$, ensuring that 
	$\Gamma$ is unbiased, i.e.,
	$\E[\Gamma] = \mu$, since
	$\E[Y_{a,b}] = \mu$ by stationarity.
	We show next that an appropriate
	choice of $\gamma_a$s produces a
	$\sqrt{n}$-consistent estimator. 
	
	Let $\mathbf{Z}_1 = (Z_{1,1},\ldots,Z_{1,\tilde{n}})$.
	By the conditional variance
	formula,
	\begin{equation}\label{eq:var-gamma}
	\var[\Gamma]
	=\var[\E[\Gamma\,|\,\mathbf{Z}_1]]
	+ \E[\var[\Gamma\,|\,\mathbf{Z}_1]].
	\end{equation}
	We first choose the $\gamma_a$s such that
	$\E[\Gamma\,|\,\mathbf{Z}_1] = \mu$;
	hence the first term on the right-hand side of~\eqref{eq:var-gamma}
	is $0$. To simplify the proof,
	we assume that the eigenvalues
	$\lambda_\ell$s are distinct (but see below). 
	By symmetry and the Markov property of the
	$(\T,P)$-walk,
	\begin{equation}
	\label{eq:gammab-1}
	\E[\Gamma\,|\,\mathbf{Z}_1]
	= \frac{1}{\tilde{n}} \sum_{b=1}^{\tilde{n}}  \E[\Gamma_b\,|\,Z_{1,b}].
	\end{equation}
	By the eigendecomposition Eq. (8) in the main text,
	we have 
	\begin{eqnarray}
	\E[\Gamma_b\,|\,Z_{1,b}]
	&=& \sum_{a=1}^{K} \gamma_a
	\sum_{\ell=1}^{N} 
	\lambda_\ell^{a-1}
	\beta_\ell
	f_\ell(Z_{1,b}) \\
	&=&\sum_{\ell=1}^{K}
	\beta_\ell
	f_\ell(Z_{1,b})
	\left\{
	\sum_{a=1}^{K} \gamma_a
	\lambda_\ell^{a-1}
	\right\}, \label{eq:gammab-2}
	\end{eqnarray}
	where we use the convention $0^0=1$.
	An appropriate choice of $\gamma_a$s can be obtained by solving the Vandermonde
	system
	\begin{equation}
	\label{eq:gammab-3}
	\sum_{a=1}^{K} \gamma_a
	\lambda_\ell^{a-1}
	= \delta_1(\ell),
	\end{equation}
	where $\delta_1(\ell) = 1$ if and only if
	$\ell = 1$, and is $0$ otherwise.
	The system above has a unique solution when the
	$\lambda_\ell$s are distinct. By \eqref{eq:gammab-1}, \eqref{eq:gammab-2}, and \eqref{eq:gammab-3}, and
	the fact that $\beta_1 = \mu$ and $f_1 = \1$, we get that
	\begin{equation}
	\label{eq:gammab-4}
	\var[\E[\Gamma\,|\,\mathbf{Z}_1]]
	= \var[\mu]
	= 0.
	\end{equation}
	
	Further, by the Markov property of the
	$(\T,P)$-walk, conditioned on $\mathbf{Z}_1$,
	the $\Gamma_b$s are independent. Hence
	$$
	\var[\Gamma\,|\,\mathbf{Z}_1]
	=\frac{1}{\tilde{n}^2} \sum_{b=1}^{\tilde{n}}  \var[\Gamma_b\,|\,Z_{1,b}],
	$$
	and, by stationarity,
	$$
	\E[\var[\Gamma\,|\,\mathbf{Z}_1]]
	= \frac{1}{\tilde{n}} \E[\var[\Gamma_b\,|\,Z_{1,b}]]
	= \frac{1}{\tilde{n}} v(P),
	$$
	where $v(P)$ is a function of $P$,
	not depending on $n$. Hence
	the second term on the right-hand side of~\eqref{eq:var-gamma}
	is $O(1/n)$. Together with \eqref{eq:gammab-4}, we have finally
	$$
	\var[\Gamma]
	=  \frac{1}{\tilde{n}} v(P)
	= O(1/n).
	$$
	That concludes the proof.
	
	If the eigenvalues are not distinct,
	let $\tilde{K}$ be the number of distinct
	eigenvalues in $\mathcal{K}$. Let $\Lambda_1,\ldots,\Lambda_{\tilde{K}}$
	be the corresponding eigenvalues. Further, if 
	$\Lambda_{i'} = \lambda_i = \cdots = \lambda_{i+k}$,
	let $g_{i'} = \beta_i f_i + \cdots + \beta_{i+k} f_{i+k}$. Then solve the reduced Vandermonde
	system so obtained.
\end{proof}

\subsection{Proof of Proposition~\ref{prop:Q}}

\begin{proof}
	Because (i) the $(\T,\sP)\textit{-walk on } G$ is stationary, (ii) $\hat Q$ can be written as a sum over steps of the chain, and (iii) the expectation is linear, it suffices to argue about one step of the chain. 
	Throughout this proof, probabilities and expectations are with respect to the $(\T,\sP)\textit{-walk on } G$.
	
	Under the DC-SBM, the expectation of the matrix $A$  has the form
	\[\A = \Theta Z{\bf B}Z^T \Theta,\]
	where $\Theta$ is a diagonal matrix whose $(i,i)$-th element is $\theta_i$. So, under the assumption that ${\bf B}$ is symmetric, so is $\A$.  By standard results for random walks on weighted graphs (see, e.g., \cite{levin2009markov}), this makes $\sP = \D^{-1} \A$ a reversible transition matrix with stationary distribution 
	\[\pi_i^* = \frac{\D_{ii}}{\sum_j \D_{jj}},\]
	where recall that $\D \in \R^{N \times N}$ is a diagonal matrix with $(i,i)$-th element $\sum_j \A_{ij}$.
	Under the stationary distribution, the probability of a transition from node $i$ to node $j$ has probability 
	\begin{equation}\label{eq:aijm}
	\pr(X_t = i, X_{t+1} = j) = \pi_i^* \pr(X_{t+1} = j| X_t = i) = \A_{ij}/m,
	\end{equation}
	where $m =\sum_j \D_{jj}$.  Using the fact that the $\theta_i$s sum to one within each block,
	\begin{equation}\label{eq:defm}
	m = \sum_j \D_{jj} = \sum_{ij} \A_{ij} = \1^T \A \1 = \1^T \Theta Z {\bf B} Z^T \Theta \1 =  \1_K^T {\bf B} \1_K.
	\end{equation}
	Again using the fact that the $\theta_i$s sum to one within each block,
	\begin{equation} \label{eq:zsum}
	(Z^T \Theta Z)_{u,v}
	=
	\sum_{i,j} Z_{i,u} \Theta_{i,j} Z_{j,v}
	=
	\sum_{i} Z_{i,u} \Theta_{i,i} Z_{i,v}
	= \sum_{i} \theta_i Z_{i,u} Z_{i,v},
	\end{equation}
	which is $1$ exactly when $u = v$ and zero otherwise.  That is, $Z^T \Theta Z$ is the identity.

	Denote $z(X_i)$ as the block membership of the $i$-th node in the Markov chain.   Note that, by definition, element $(u,v)$ of $Q = \E(\hat Q)  \in \R^{K \times K}$ is the probability of a transition from block $u$ to block $v$. For $S(u,v) = \{(i,j): z(i) = u, z(j) = v\}$,
	\begin{eqnarray*}
		Q_{u,v}&=&\pr\left(z(X_t) = u, z(X_{t+1}) = v\right)  \\
		&=& \sum_{ (i,j) \in S(u,v)} \pr(X_t = i, X_{t+1} = j) \\
		&=& \sum_{(i,j) \in S(u,v)} \A_{ij}/m,
	\end{eqnarray*}
	by \eqref{eq:aijm}.
	%
	With \eqref{eq:zsum}, this implies the conclusion,
	\[Q  = (1/m) Z^T \A Z = (1/m) (Z^T  \Theta Z){\bf B}(Z^T \Theta Z) = (1/m) {\bf B}.\]

\end{proof}

\subsection{Proof of Proposition~\ref{prop:B}}

\begin{proof}

	Because ${\bf B}$ is symmetric, $B_L = D_B^{-1/2} {\bf B}D_B^{-1/2}$  (with $D_B = \diag( {\bf B} \1_K)$) has a real-valued eigendecomposition, $B_L = U \Lambda U^T$, where $U \in \R^{K \times K}$ is orthonormal and $\Lambda \in \R^{K \times K}$ is diagonal.  By Lemma 3.3 in \cite{qin2013regularized}, the eigenvalues of $B_L$ match the nonzero eigenvalues of $\sP$.  Moreover, this lemma shows that 
	the orthonormal columns of $U^* = \Theta^{1/2}ZU \in \R^{N \times K}$  are the  eigenvectors of $\D^{-1/2}\A \D^{-1/2}$.  Left multiply both sides of the eigenvector equations  by $\D^{-1/2}$:
	\[\D^{-1/2} [ \D^{-1/2}\A \D^{-1/2} U^*]=
	\D^{-1/2}[U^* \Lambda] \implies\]
	\[\D^{-1}\A [\D^{-1/2} U^*]=
	[\D^{-1/2}U^*] \Lambda.
	\]
	That is, the columns of $\D^{-1/2} U^*$ are eigenvectors of $\sP = \D^{-1} \A$. Because $\D_{ii} =  \theta_i [D_B]_{z(i), z(i)}$ and the relationship between $z(i)$ and $Z$,
	\[[\D^{-1/2} \Theta^{1/2}ZU]_{ij}  = [D_B]^{-1/2}_{z(i), z(i)} [ ZU]_{ij} 
	= [ZD_B^{-1/2}U]_{ij}.\]
	Thus, the eigenvectors of $\sP$ simplify to $Z D_B^{-1/2} U$. We must still adjust the scaling  to ensure that the eigenvectors $f_v^*$ satisfy, 
	\begin{equation}\label{eq:normalization} 
	\sum_{i = 1}^N f_v^*(i) f_u^*(i) \pi_i^* = \left\{\begin{array}{cc} 1 & u=v \\0 & o.w.\end{array}\right.,
	\end{equation}
	where $\pi^*$ is the stationary distribution of $\sP$.  Because $\A$ is symmetric, the chain is reversible and $\pi^*$ contains the normalized row sums of $\A$, \[\pi_i^* = \frac{\D_{ii}}{\sum_j D_{jj}} = \frac{\theta_i [D_B]_{z(i), z(i)} }{\sum_j D_{jj}}.\]
	Denote $\Pi^* \in \R^{N \times N}$ as the diagonal matrix that contains $\pi^*$ down the diagonal.  
	
	Define $f^* =\sqrt{m} Z D_B^{-1/2} U \in \R^{N \times K},$ where $m = \1_K^T {\bf B} \1_K$ and define $f_v^*$ as the $v$th column of $f^*$.  Note that the $u,v$th element of the following matrix is the sum defined on the left side of the \eqref{eq:normalization},
	\begin{eqnarray}
	\label{eq:fpf} [f^*]^T \Pi^* f^* &=& [\sqrt{m} Z D_B^{-1/2} U]^T \Pi^* \sqrt{m} Z D_B^{-1/2} U \\
	&=& m \ U^T  D_B^{-1/2} Z^T \Pi^* Z D_B^{-1/2} U. \label{eq:fpf2} 
	\end{eqnarray}
	To simplify this term, recall that the definition of the Degree Corrected Stochastic Blockmodel makes the model identifiable by presuming
	\[\sum_{i:z(i) = u} \theta_i = 1\]
	for all blocks $u$.  This implies that 
	\[Z^T \Pi^* Z = D_B (\sum_j D_{jj})^{-1} = m^{-1} \ D_B .\] Thus, \eqref{eq:fpf2} simplifies to the identity matrix. 
	
	The last piece of the Theorem follows from definitions,
	\[ \beta_\ell^* = \langle y, f_\ell^* \rangle_\pi = \sum_i y(i) f_\ell^*(i) \pi_i  = \E(y(X) f_\ell^*(X)). \]
	
\end{proof}

\section{A diagnostic plot}   \label{sec:diagnostic}
To examine and compare the different fGLS estimators, this main text proposes a diagnostic plot.  Importantly, the plot does not require any additional information beyond the information used to compute the fGLS estimators.  As such, it can be used in practice to help compare different fGLS estimators.  The plot is explained again below (in slightly more detail).  Then, there is some theoretical support for these diagnostics given in the rest of the section.

\begin{figure}[h] 
	\centering
	\includegraphics[width=5in]{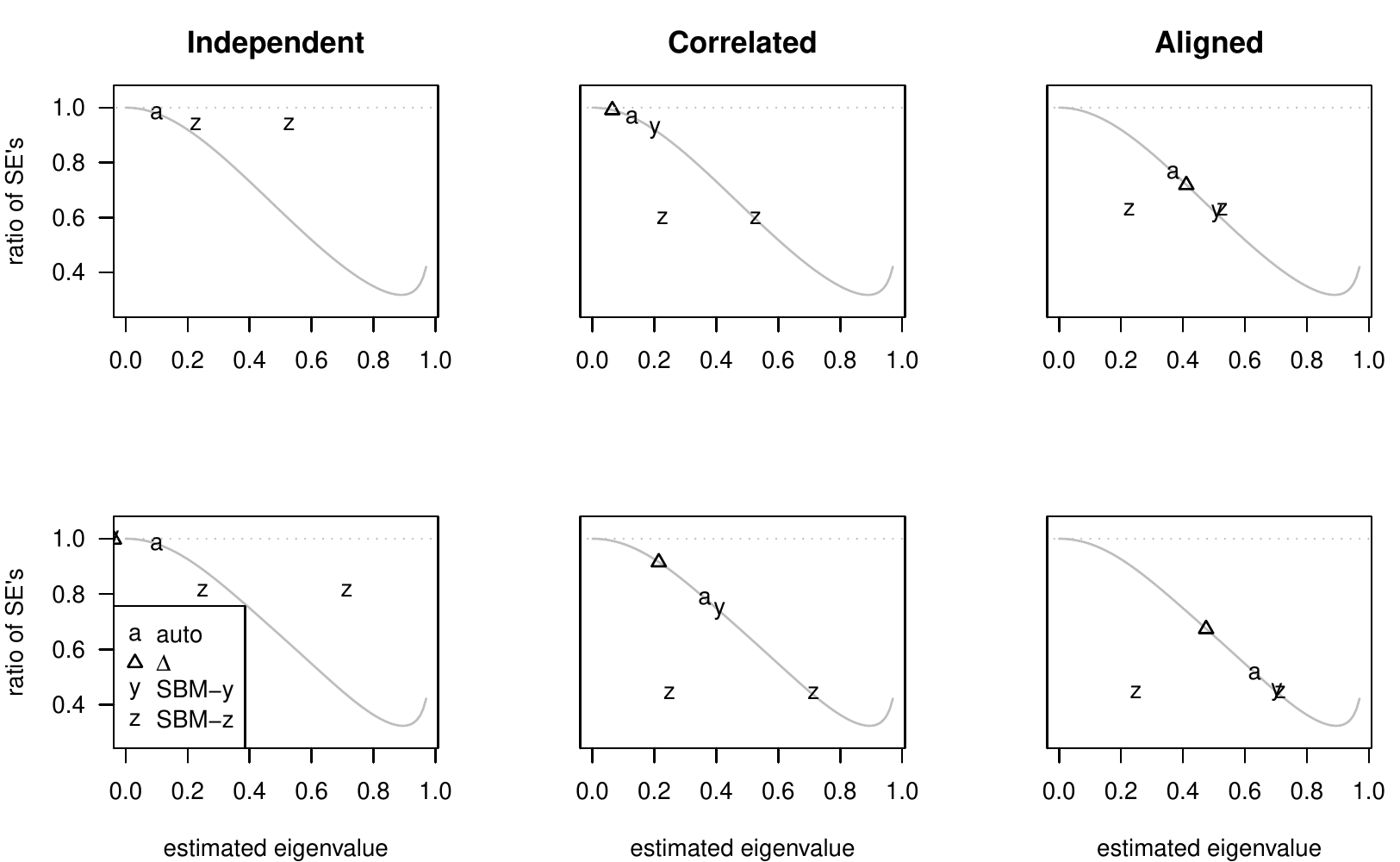} 
	\caption{
		Each of these diagnostic plots is created from a single sample of $n=500$ from the respective panel in Figure \ref{fig:rmse}.  Based upon the presumptions described and justified in the rest of this section, we should prefer estimators that are further to the right and further to the bottom because they detect more dependencies and thus make further reductions in the variance.  For example, the middle panel in the bottom row shows a clear preference for the SBM-z estimator.  Because these are simulations, we know that this is the only estimator that correctly specifies the model; consulting the Figure \ref{fig:rmse} this is the estimator that has the smallest RMSE.  
	}
	\label{fig:diagnostics-si}
\end{figure}

In the diagnostic plot, each rank-two fGLS estimator appears as a single point and each SBM-fGLS estimator appears as $K-1$ points, where $K$ is the number of blocks in the SBM.  
The horizontal axis of the plot gives estimated eigenvalue(s) of $P$ and the vertical axis gives the estimated \textbf{ratio of standard errors},  
\begin{equation}\label{eq:rse}
\rse(\hat \Sigma) = 
\sqrt{
	\frac{\widehat{\var}(\hat \mu_{\gls}(\hat \Sigma))} 
	{\widehat{\var}(\hat \mu)}
} = 
\sqrt{\frac
	{(\1^T \hat \Sigma^{-1} \1)^{-1} }
	{n^{-1} \1^T \hat \Sigma \1}},
\end{equation}
where the numerator and denominator are plug-in estimates of the standard errors, both using the same estimated covariance matrix $\hat \Sigma$.
Because different estimators construct different estimates $\hat \Sigma$, each estimator has a different $\rse$, represented in the diagnostic plots with the symbols {\bf a, $\Delta$, y} and ${\bf z}$.  In the plot, there are $(K-1)$-many ${\bf z}$s because SBM-fGLS estimates $K-1$ eigenvalues; each of these $K-1$ points has the same value on the vertical axis.  Each of the six panels in Figure \ref{fig:diagnostics-si} gives the diagnostic plot for the first sample of $n=500$ that appeared in the respective panel of Figure \ref{fig:rmse-si}.  The grey line in the diagnostic plot is described below.

It is important to note that the $\rse$ is likely to be a biased estimate of the actual reduction in standard error because the plug-in estimator for $\var(\hat \mu_{\gls}(\hat \Sigma))$ does not account for the uncertainty in estimating $\hat \Sigma$.  Moreover, SBM estimators with a large $K$ are likely to overfit, leading to misleadingly low values of $\rse$.  

\subsection{Interpreting the diagnostic plot relies on two key presumptions}
The following  two presumptions motivate the diagnostic plots: 
\begin{enumerate}
	\item[\textbf{P1:}]  If an estimator is based upon a misspecified model for $\Sigma$, then it will under specify the dependence between samples.  As such, the estimated eigenvalue(s) and $n^{-1} \1^T \hat \Sigma \1$ will be smaller than in a correctly specified model.
	\item[\textbf{P2:}] If a model does not detect the full dependence in $\Sigma$, then the fGLS estimator will not reduce the variance as much as the fGLS based upon the correctly specified model.  As such, the $\rse$ will be larger than in a correctly specified model. 
\end{enumerate}

\subsubsection{The key evidence for \textbf{P1} comes from a result in  \cite{li2015central}, restated below in the notation of the current paper} 

\ Define 
\[\lambda_{\auto} = \frac{\gamma(1)}{ \gamma(0)}\] as the population version of $\hat \lambda$ in \eqref{eq:lamauto}.  Under the rank-two model, $\lambda_{\auto} = \lambda_2(P)$.  When the rank-two model is misspecified, $\hat \mu_{\auto}$ still relies on estimating $\lambda_{\auto}$, which is not necessarily equal to $\lambda_2(P)$.  
Define  $\Sigma^{(\auto)} \in \R^{n \times  n}$ as  indexed by the nodes in $\T$,
\begin{equation} \label{eq:SigmaAuto}
\Sigma^{(\auto)}_{\sigma, \tau} = \gamma(0) \lambda_{\auto}^{d(\sigma,\tau)},
\end{equation}
where $d(\sigma,\tau)$ is the distance between  $\sigma$ and $\tau$ in $\T$.  

\begin{prop}[From \cite{li2015central}] \label{prop:jensen}
	If $P$ is reversible and irreducible then $|\lambda_{\auto}| \le |\lambda_2(P)|$.  If it is also true that all eigenvalues of $P$ are non-negative,  then 
	\begin{equation}\label{eq:jensen}
	\1^T \Sigma \1 \ge \1^T \Sigma^{(\auto)}\1.
	\end{equation}
	Equality holds when $\gamma$ is rank-two.
\end{prop}
The proof is based on Jensen's inequality.  For completeness,  the proof is included at the end of this section.  See \cite{li2015central} for a discussion of how the argument often extends to allow $P$ to have negative eigenvalues.

The left side of \eqref{eq:jensen} divided by $n$ is the variance of the sample average.\footnote{As discussed in Remark 1.1, this argument extends to unbiased estimates as well.} Under the rank-two model, \eqref{eq:jensen} holds with equality. As such, $\Sigma^{(\auto)}$ provides one way to estimate the variance; the denominator of $\rse$ for $\hat \mu_{\auto}$ relies upon this fact.  
As such, when the rank-two model is misspecified, 
$|\lambda_{\auto}| \le |\lambda_2(P)|$ and \eqref{eq:jensen} imply that $\Sigma^{(\auto)}$ understates the dependence. This supports \textbf{P1}.  



%
%
%
%
%

\subsubsection{The key evidence for \textbf{P2} comes from the grey line in the diagnostic plot} \ 

Let $\tilde \Sigma$ be a covariance matrix based upon the rank-two model for the autocovariance function. That is, for some $\tilde \beta$ and $\tilde \lambda$,
\begin{equation} \label{eq:rankTwoSigma}
\tilde \Sigma_{\sigma, \tau} = \tilde \beta \tilde \lambda^{d(\sigma,\tau)}.
\end{equation}
The argument below shows that $\rse(\tilde \Sigma)$ is completely determined by $\tilde \lambda$ and the structure of $\T$, which is constant across estimators.  The diagnostic plot uses this fact by giving $\rse(\tilde \Sigma)$ as a function of $\tilde \lambda$.  The key evidence for \textbf{P2} is that in Figure \ref{fig:diagnostics} the grey line is downward sloping in the region $\tilde \lambda \in [0, .9]$.  As such, estimating $\lambda_2(P)$ closer to zero will yield a larger $\rse$.

To see that $\rse(\tilde \Sigma)$ does not rely on $\tilde \beta$, first look at the numerator of $\rse$.   Rearrange the terms in \eqref{eq:1sig1} to see that for some function $s_{gls}$,
\[(\1^T \tilde \Sigma^{-1} \1)^{-1} = \tilde\beta^{2} s_{gls}(\tilde \lambda, n).\]
Interestingly, $s_{gls}$ does not depend on the topological properties of  $\T$, only the number of observations $n$. 
Next, look at the denominator of $\rse$.  For some function $s$,
\[ n^{-1}\1^T\tilde \Sigma \1 = n^{-1}\sum_{\tau, \sigma} \tilde \beta^2 \tilde  \lambda^{d(\tau,\sigma)} = \tilde  \beta^2 s(\tilde \lambda, \T).\]
Thus, under the rank-two model,
\[\rse(\tilde  \Sigma)^2  = 
{\frac
	{(\1^T \tilde  \Sigma^{-1} \1)^{-1} }
	{n^{-1} \1^T \tilde  \Sigma \1}} 
=
{\frac{\tilde \beta^{2} s_{gls}(\tilde \lambda,n)}{\tilde \beta^2 s(\tilde \lambda, \T)}}
=
{\frac{s_{gls}(\tilde \lambda,n)}{s(\tilde \lambda, \T)}}.\]
Thus, for any fixed $\T$, $\rse$ is only a function of $\tilde \lambda$.

%

\subsection{Proof of Proposition \ref{prop:jensen}}
\begin{proof}
	First, to see that $|\lambda_{\auto}| \le |\lambda_2(P)|$, define $w_\ell = \beta_\ell^2 / \sum_{j=2}^N \beta_j^2$.  Because $|\lambda_2(P)| \ge |\lambda_\ell(P)|$ for all $\ell\ge 2$, 
	\[|\lambda_{\auto}| = \left|\frac{\gamma(1)}{\gamma(0)}\right| =  \left|\sum_{\ell = 2}^N w_\ell \lambda_\ell(P)  \right| \le  \sum_{\ell = 2}^N w_\ell |\lambda_\ell(P)| \le |\lambda_2(P)|.\]

	By Theorem 2.1 in  \cite{treevar}, 
	\[n^{-1}\1^T \Sigma \1  = \sum_{\ell = 2}^N \beta_\ell^2 \G(\lambda_\ell),\]
	where $\G$ is defined as follows.  Let $I$ and $J$ be  drawn independently and uniformly at random from the nodes of $\T$.  Define the random variable $D = d(I,J)$ as the (random) distance.  Define $\G$ as the probability generating function for $D$,
	\[\G(x) = \E(x^D).\]
	Because $\G$ is the probability generating function of a non-negative and discrete random variable, $\G$ is convex on $[0,1)$ (e.g. p29 in \cite{durrett2007random}).  By assumption $\lambda_\ell \ge 0$ for $\ell \in 1, \dots, N$.  Moreover,  because $P$ is irreducible,  $\lambda_2  <1$.  So, by Jensen's inequality,
	\begin{eqnarray*}
		n^{-1}\1^T \Sigma \1 
		&=& \sum_{\ell = 2}^N \beta_\ell^2 \G(\lambda_\ell)\\
		&=& \left(\sum_{j=2}^N \beta_j^2\right) \sum_{\ell = 2}^N \frac{\beta_\ell^2}{\sum_{j=2}^N \beta_j^2} \G(\lambda_\ell)\\
		&\ge& \left(\sum_{j=2}^N \beta_j^2 \right) \G\left(\sum_{\ell = 2}^N \frac{\beta_\ell^2 \lambda_\ell}{\sum_{j=2}^N \beta_j^2}\right)\\
		&=& \gamma(0)  \ \G\left(\frac{\gamma(1)}{\gamma(0)}\right)\\
		&=& n^{-1}\1^T \Sigma^{(\auto)}\1,
	\end{eqnarray*}
	where the last line again follows from Theorem  2.1 in \cite{treevar}.
\end{proof}

\section{Construction of estimators with SBM-fGLS estimation of normalization constant} \label{sec:construction}

This section describes the construction of the estimators in Section~\ref{sec:fgls} in the main paper.
For notational convenience, denote $Y_\tau$, $\tilde z(\tau)$, and $\deg(\tau)$ as $y(X_\tau), \tilde z(X_\tau),$ and $\deg(X_\tau)$ for each sampled individual $X_\tau$. Moreover, suppose a one-to-one mapping  between the node set of $\T$ and $\{1, \dots, n\}$.

The key function $\texttt{sbm-fgls}(\T, \{Y_\tau, \tilde z(\tau)\}_{\tau \in \T})$ is defined in the next subsection.  It returns an SBM-fGLS point estimate for $\E(Y_\tau)$.

The point estimates in Section~\ref{sec:fgls} are computed as follows. 
\begin{enumerate}
	\item[1)] Estimate $\E(1/\deg(X_\tau))$ with $H_{gls}^{-1} = \texttt{sbm-fgls}(\T, \{1/\deg(\tau), \tilde z(\tau)\}_{\tau \in \T})$.  
	\item[2)] Compute 
	\[Y_\tau^{\hat \pi} = \frac{Y_\tau}{H_{gls}^{-1} \deg(\tau)}.\]
	\item[3)] Return $\texttt{sbm-fgls}(\T, \{Y_\tau^{\hat \pi}, \tilde z(\tau)\}_{\tau \in \T})$.
\end{enumerate}

\subsection{The definition of \texttt{sbm-fgls}}
To compute $\texttt{sbm-fgls}(\T, \{Y_\tau, \tilde z(\tau)\}_{\tau \in \T})$, follow these steps:

\begin{enumerate}
	\item[1)] Compute $\hat Q$ via Equation (14) in the main text using the block memberships $\tilde z(\tau)$.  Define $\hat Q^{(S)} = (\hat Q + \hat Q^T)/2$. This symmetrization ensures the  eigenvalues are real-valued. 
	\item[2)] Row and column normalize $\hat Q^{(S)}$, as $\hat Q_L = D_{\hat Q}^{-1/2} \hat Q^{(S)} D_{\hat Q}^{-1/2},$
	where $D_{\hat Q} = \diag(\hat Q \1_K) \in \R^{K \times K}$.
	\item[3)] Take an eigendecomposition of 
	\begin{equation}\label{eq:lamsbm}
	\hat Q_L = \hat U \hat \Lambda \hat U^T.
	\end{equation}
	\item[4)] Compute $\hat f = \hat Z D_{\hat Q}^{-1/2} \hat U$, where $\hat Z \in \{0,1\}^{n\times K}$ contains $\hat Z_{ij} = 1$ iff $\tilde z(i) = j$. 
	\item[5)] For $\ell = 1, \dots, K$, compute $\hat \beta_\ell = \frac{1}{n} \sum_\tau Y_\tau \hat f_\ell(\tau)$, where $\hat f_\ell(\tau)$ is the $(\ell, \tau)$ element of $\hat f$.
	\item[6)] Compute an estimate of the auto-covariance function as
	\[\hat \gamma_{\sbm}(d) = \sum_{\ell = 1}^K  \hat \beta_\ell^2 \hat \Lambda_{\ell \ell}^d.\]
	\item[7)] Define $\hat s^2$ to be the sample variance of the $Y_\tau$.  
	For $\sigma, \tau \in \T$,
	$$\hat \Sigma_{\sigma, \tau}^{\sbm} = 
	\left\{\begin{array}{lc}
	\hat \gamma_{\sbm}(d(\sigma, \tau)) & \text{if $\sigma \ne \tau$} \\
	\hat \gamma_{\sbm}(0) + \hat s^2& \text{if $\sigma  = \tau$,} \end{array}\right.$$
	where $\hat s^2$ provides for Tikhonov regularization in $(\hat \Sigma^{\sbm})^{-1}$.
	\item[8)] Define $\hat g \in \R^n$ to solve the system of equations $\hat \Sigma_{sbm} \hat g = \textbf{1}$.
	\item[9)] Estimate $\E(Y_\tau)$ with $\sum_{\tau \in \T} \hat g_\tau Y_\tau / \sum_{\tau \in \T} \hat g_\tau$.\end{enumerate}

\end{document}